
\input amstex
\documentstyle{amsppt}

\pagewidth{6.5truein}
\pageheight{9.375truein}
\vcorrection{-0.375truein}

\NoBlackBoxes
\TagsAsMath
\TagsOnRight

\long\def\ignore#1\endignore{#1}

\ignore
\input xy
\xyoption{matrix}\xyoption{arrow}\xyoption{curve}\xyoption{frame}
\def\edge{\ar@{-}}

\def\wiggledge{\ar@{~}}
\def\dropup#1#2{\save+<0ex,#1ex> \drop{#2} \restore}
\def\dropdown#1#2{\save+<0ex,-#1ex> \drop{#2} \restore}
\def\encirc#1{*+[o][F]{#1}}
\endignore

\def\la{\Lambda}
\def\lamod{\Lambda{}\operatorname{-mod}}
\def\deltamod{\Delta{}\operatorname{-mod}}
\def\Deltamod{\Delta{}\operatorname{-Mod}}
\def\Lamod{\Lambda{}\operatorname{-Mod}}

\def\pinflamod{{\Cal P}^{< \infty}(\Lambda{}\operatorname{-mod})}
\def\Pinflamod{{\Cal P}^{< \infty}(\Lambda{}\operatorname{-Mod})}
\def\pinfdeltamod{{\Cal P}^{< \infty}(\Delta{}\operatorname{-mod})}
\def\Pinfdeltamod{{\Cal P}^{< \infty}(\Delta{}\operatorname{-Mod})}

\def\sinflamod{{\Cal S}^{< \infty}(\Lambda{}\operatorname{-mod})}
\def\sinfdeltamod{{\Cal S}^{< \infty}(\Delta{}\operatorname{-mod})}

\def\pdim{\operatorname{p\,dim}}
\def\Hom{\operatorname{Hom}}

\def\lfindim{\operatorname{l\,fin\,dim}}

\def\lFindim{\operatorname{l\,Fin\,dim}}

\def\Im{\operatorname{Im}}
\def\length{\operatorname{length}}
\def\soc{\operatorname{soc}}
\def\supp{\operatorname{supp}}
\def\NN{\Bbb N}
\def\ZZ{\Bbb Z}

\def\A{{\Cal A}}
\def\C{{\Cal C}}
\def\D{{\Cal D}}
\def\Add{\operatorname{Add}}

\def\ehat{\hat{e}}

\def\phat{\hat{p}}
\def\qhat{\hat{q}}
\def\vhat{\hat{v}}
\def\hatp{\hat{p}}
\def\hatq{\hat{q}}
\def\what{\hat{w}}
\def\etil{\tilde{e}}
\def\ptil{\tilde{p}}
\def\qtil{\tilde{q}}
\def\wtil{\tilde{w}}
\def\fp{\frak{p}}
\def\fq{\frak{q}}
\def\St{\operatorname{St}}
\def\Bd{\operatorname{Bd}}

\def\AuRe{{\bf 1}}
\def\AuSm{{\bf 2}}
\def\AuSmtwo{{\bf 3}}
\def\BenGna{{\bf 4}}
\def\Bleone{{\bf 5}}
\def\Bon{{\bf 6}}
\def\BrSk{{\bf 7}}
\def\BuRi{{\bf 8}}
\def\ChriS{{\bf 9}}
\def\CraBoe{{\bf 10}}
\def\CBser{{\bf 11}}
\def\DoFr{{\bf 12}}
\def\Erd{{\bf 13}}
\def\ErSk{{\bf 14}}
\def\Gabriel{{\bf 15}}
\def\Geistring{{\bf 16}}
\def\GePo{{\bf 17}}
\def\HaHZ{{\bf 18}}
\def\pre{{\bf 19}}
\def\dom{{\bf 20}}
\def\dep{{\bf 21}}
\def\Huiathens{{\bf 22}}
\def\HZSmdnsn{{\bf 23}}
\def\HuSm{{\bf 24}}
\def\IgSmTo{{\bf 25}}
\def\Kratwo{{\bf 26}}
\def\Kra{{\bf 27}}
\def\Len{{\bf 28}}
\def\Rin{{\bf 29}}
\def\Rintwo{{\bf 30}}
\def\Ringen{{\bf 31}}
\def\RinSchroe{{\bf 32}}
\def\Rog{{\bf 33}}
\def\Schroe{{\bf 34}}
\def\SkWa{{\bf 35}}
\def\Sma{{\bf 36}}
\def\Sze{{\bf 37}}
\def\VFCB{{\bf 38}}
\def\WaWa{{\bf 39}}

\topmatter

\title The homology of string algebras I
\endtitle

\author B. Huisgen-Zimmermann  and S. O. Smal\o \endauthor

\address Department of Mathematics, University of California, Santa Barbara,
CA 93106, USA\endaddress
\email birge\@math.ucsb.edu\endemail

\address Department of Mathematics, The Norwegian University for Science and
Technology, 7055 Dragvoll, Norway
\endaddress
\dedicatory Dedikert til v{\aa}r venn og kollega Idun Reiten i anledning
hennes seksti{\aa}rsdag \enddedicatory
\email sverresm\@math.ntnu.no \endemail

\dedicatory Dedikert til v{\aa}r venn og kollega Idun Reiten i andledning
hennes seksti{\aa}rsdag \enddedicatory
\thanks The research of the first author was partially supported by a grant
from the National Science Foundation.  Both authors would like to
acknowledge the hospitality of the Center for Advanced Study of the
Norwegian Academy of Sciences and Letters where the project was completed.
\endthanks

\abstract  We show that string algebras are `homologically tame' in the
following sense:  First, the syzygies of arbitrary representations of a
finite dimensional string algebra $\la$ are direct sums of cyclic
representations, and the left finitistic dimensions, both little and big, of
$\la$ can be computed from a finite set of cyclic left ideals contained in
the Jacobson radical.  Second, our main result shows that the functorial
finiteness status of the full subcategory
$\pinflamod$ consisting of the finitely generated left $\la$-modules of
finite projective dimension is completely determined by a finite number of,
possibly infinite dimensional, string modules  --  one for each simple
$\la$-module  --  which are algorithmically constructible from  quiver and
relations of $\la$.  Namely, $\pinflamod$ is contravariantly finite in
$\lamod$ precisely when all of these string modules are finite dimensional,
in which case they coincide with the minimal
$\pinflamod$-approximations of the corresponding simple modules.  Yet, even
when
$\pinflamod$ fails to be contravariantly finite, these `characteristic'
string modules encode, in an accessible format, all desirable homological
information about $\lamod$.
\endabstract
\endtopmatter

\document

\head 1. Introduction \endhead

The representation theory of the Lorentz group is intimately linked to that
of a certain string algebra $\la$ (for a definition of string algebras see
Section 2), as was observed and exploited by Gelfand and Ponomarev in
\cite{\GePo}.  In particular, it was proved there that this algebra  --
along with a class of close relatives  --  has tame representation type; in
fact, its finite dimensional indecomposable representations were explicitly
pinned down.  In a sequence of articles by Ringel \cite{\Rin}, Bondarenko
\cite{\Bon}, Donovan-Freislich \cite{\DoFr}, Butler-Ringel
\cite{\BuRi}, and others, the class of algebras amenable to techniques
derived from the Gelfand-Ponomarev archetype was subsequently found to be
much larger and, moreover, to be related to further classical scenarios,
such as the representation theory of dihedral groups.  This development
ultimately led to a well-rounded representation-theoretic picture of the
extended class of algebras on which we concentrate here, the class of string
algebras.  Among other tools, Auslander-Reiten methods were employed to place
the finite dimensional indecomposable objects into a tightly knit
categorical context.  In tandem, certain portions of the infinite
dimensional representation theory were rendered accessible. However, in
spite of the availability  of a full classification of the finitely
generated indecomposable representations of string algebras, their
homological properties, known to vary widely (see, e.g., \cite{\IgSmTo}),
remained far from understood.

Our goal is to supplement the structural information with equally precise
homological data.  For a more detailed preview of our results, let $\la$ be
a finite dimensional string algebra over an algebraically closed field $K$. 
We start by showing how the homological dimensions --  global and
finitistic  --  can be obtained from a finite collection of cyclic modules
contained in the Jacobson radical $J$ of $\la$ (Theorem 3).  Then we turn to
two far more deep-seated problems concerning the category 
$\pinflamod$ that has as objects the finitely generated modules
of finite projective dimension.  These problems are as follows: (I)
Can the internal structure of the objects in $\pinflamod$ be characterized, so as
to distinguish them from those in $\lamod
\setminus \pinflamod$?  --  here $\lamod$ stands for the category of all
finitely generated left $\la$-modules  --  and (II) How is the category
$\pinflamod$ embedded in
$\lamod$, in terms of maps entering or leaving $\pinflamod$? It is our
answer to the second question which displays the homological mechanisms of
$\la$; in particular, it entails a solution to the first problem.

To tackle Problem II, we establish a readily checkable characterization of
contravariant finiteness of the category
$\pinflamod$ in $\lamod$ and describe the minimal
$\pinflamod$-approximations of the simple modules $S_1, \dots, S_n$ in the
positive case (cf\. Section 4 for the relevant definitions); in fact, our
description yields a procedure for constructing them (Theorem 5 and
Proposition 14). Existence of such minimal $\pinflamod$-approximations of
the $S_i$ and  --  existence provided  --  their structure are known to have
far-reaching consequences for the homology of $\la$ (see, e.g.
\cite{\AuSmtwo},
\cite{\AuRe}, \cite{\HuSm}); those which have direct impact on our
present investigation are reviewed in Section 4.

To appreciate how contravariant finiteness of $\pinflamod$ relates to our
second problem, recall that this condition implies functorial finiteness of
$\pinflamod$, i\.e., dually defined left
$\pinflamod$-approximations of the objects in $\lamod$ come as free
byproducts (see \cite{\HuSm}). Suppose, for the moment, that
$\pinflamod$ is contravariantly finite, and fix an object $M \in
\lamod$.  The key roles played by the minimal right
$\pinflamod$-approximation $\phi_M: A_M \rightarrow M$ and the minimal left
approximation $\psi_M: M \rightarrow B_M$ of $M$ can be visualized as
follows:

$$
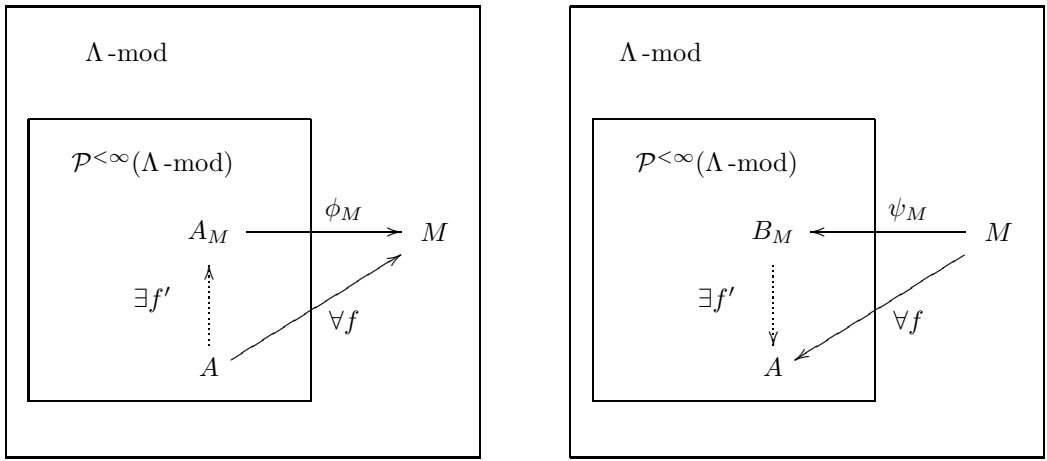
\xy{ {(16,54)*={\lamod}},{{(0,60)}="+";{(63,60)}="*" **@{-}},
{{(63,60)}="+";{(63,0)}="*" **@{-}}, {{(0,0)}="+";{(63,0)}="*" **@{-}},
{{(0,0)}="+";{(0,60)}="*" **@{-}}, {{(3,7.5)}="+";{(3,45)}="*" **@{-}},
{{(40.5,45)}="+";{(3,45)}="*" **@{-}}, {{(40.5,45)}="+";{(40.5,7.5)}="*"
**@{-}},{{(3,7.5)}="+";{(40.5,7.5)}="*" **@{-}},(19.5,39)*={\pinflamod},
(27,30)*={A_M}, (27,12)*={A}, {(27,15)}="+";{(27,25.5)}="*"
**@{.}?>*\dir{>},{(30,13)}="+";{(52.5,27)}="*"
**@{-}?>*\dir{>},{(32,30)}="+";{(52.5,30)}="*" **@{-}?>*\dir{>},
(57,30)*={M}, (19.5,21)*={\exists f'},(45,33)*={\phi_M},(45,18)*={\forall f},
 {(87,54)*={\lamod}},{{(75,60)}="+";{(138,60)}="*" **@{-}},
{{(138,60)}="+";{(138,0)}="*" **@{-}}, {{(75,0)}="+";{(138,0)}="*" **@{-}},
{{(75,0)}="+";{(75,60)}="*" **@{-}}, {{(78,7.5)}="+";{(78,45)}="*" **@{-}},
{{(115.5,45)}="+";{(78,45)}="*" **@{-}}, {{(115.5,45)}="+";{(115.5,7.5)}="*"
**@{-}},{{(78,7.5)}="+";{(115.5,7.5)}="*" **@{-}},(94.5,39)*={\pinflamod},
(102,30)*={B_M}, (102,12)*={A}, {(102,15)}="+";{(102,25.5)}="*"
**@{.}?<*\dir{<},{(105,13)}="+";{(127.5,27)}="*"
**@{-}?<*\dir{<},{(107,30)}="+";{(127.5,30)}="*" **@{-}?<*\dir{<},
(132,30)*={M}, (94.5,21)*={\exists f'},(120,33)*={\psi_M},(120,18)*={\forall
f}}\endxy
$$

\noindent In other words, $A_M$ is minimal with the property that all
homomorphisms from any object $A \in \pinflamod$ to $M$ pass through the way
station $A_M$ via $\phi_M$, and consequently, the problem of controlling all
maps in $\Hom_{\la}(\pinflamod, M)$ boils down to describing the
approximation $\phi_M$ and understanding the internal maps of $\pinflamod$.
The map $\psi_M$ plays a dual role.  On the side, we mention that functorial
finiteness of $\pinflamod$ guarantees the existence of almost split
sequences in that category (see \cite{\AuSm} and
\cite{\AuSmtwo}).

Whether or not $\pinflamod$ is contravariantly finite in $\lamod$, one can
associate with each simple $\la$-module $S_i$ a representation $A_i$ in the
category $\Pinflamod$ of all, not necessarily finite dimensional, left
$\la$-modules of finite projective dimension, together with a canonical map
$\varphi_i: A_i \rightarrow S_i$, which is indicative of the map-theoretic
`location' of $S_i$ relative to $\pinflamod$ (Theorem 5).   In fact,
$\pinflamod$ is contravariantly finite if and only if all of the
$A_i$ are finite dimensional, in which case the maps $\varphi_i$ are the
minimal $\pinflamod$-approximations of the simple modules. Otherwise, the
$A_i$ are still `phantoms' of $\pinflamod$-approximations of the $S_i$ in the
sense of
\cite{\HaHZ} (see Section 4 for a reminder). Roughly, this means that each
$A_i$ exhibits, in the tightest possible format, the relations of those
modules $M \in \pinflamod$ which map onto $S_i$; for instance, it is
precisely when $S_i$ has finite projective dimension that $A_i \cong S_i$. 
All of these constructions are algorithmic.  Indeed, even when the modules
$A_i$ are infinite dimensional, they can be constructed in a predictable
number of steps, growing polynomially with $\dim_K
\la$, due to their periodicity properties.

In the present work, we primarily focus on the homology of the `little'
category $\lamod$  --  even though we need to resort to infinite dimensional
modules to fully understand the latter  -- whereas part II
will address additional phenomena arising in the homology of the `big' category
$\Lamod$.

As mentioned at the outset, the class of string algebras has developed into
a showcase for representation-theoretic methods, thus attesting to the
`state of the art' on various fronts.  Moreover, algebras degenerating to
string algebras play a pivotal role in understanding more general classes of
algebras.  To tie the present investigation into the context of
existing work, we add a relatively short, chronologically ordered list of further
references providing historical background and samples of different lines of
approach:
\cite{\Sze}, \cite{\Gabriel}, \cite{\SkWa}, \cite{\WaWa},
\cite{\Erd}, \cite{\CraBoe}, \cite{\CBser}, \cite{\Rintwo},
\cite{\Bleone}, \cite{\Geistring}, \cite{\Kratwo}, \cite{\Rog},
\cite{\VFCB},
 \cite{\Schroe}, \cite{\Ringen}, \cite{\RinSchroe},  \cite{\BrSk}.

\head 2. Prerequisites and conventions \endhead

Consistently, $\la = K\Gamma / I$ will denote a finite dimensional {\it
string algebra\/} over an algebraically closed field $K$. This means that
$\la$ is a {\it monomial\/} relation algebra (that is, the admissible ideal $I$ of
the path algebra
$K\Gamma$ can be generated by paths), and $\la$ is
{\it special biserial\/}.  The latter amounts to the combination of the
following two conditions on
$\Gamma$ and $I$: at most two arrows enter and at most two arrows leave any
given vertex of
$\Gamma$, and, for any arrow $\alpha$ in $\Gamma$, there is at most one arrow
$\beta$ with $\alpha \beta \notin I$ and at most one arrow $\gamma$ with
$\gamma
\alpha \notin I$.

Our convention for multiplying paths is as follows:  if
$p$ and $q$ are paths in $K\Gamma$, then $pq$ stands for `$p$ after $q$'.
Correspondingly, a {\it right subpath\/} of a path $p$ is a path $u$ such
that $p = u'u$ for some path $u'$; {\it left subpaths\/} of $p$ are defined
symmetrically.  The set of vertices of
$\Gamma$ will be identified with a full set of primitive idempotents $e_1,
\dots, e_n$ of
$\la$, and we will be referring to idempotents from this set whenever we
mention primitive idempotents.  Given any left
$\la$-module $M$, a {\it top element\/} of $M$ is an element $x$ with the
property that $ex = x$ for some primitive idempotent $e$, in which case we
will also call $x$ a top element {\it of type $e$\/} of $M$.

Throughout, $\lamod$ and $\Lamod$ will stand for the categories of finite
dimensional and arbitrary left $\la$-modules, respectively,
$\pinflamod$ will denote the full subcategory of $\lamod$ having as objects
the modules of finite projective dimension, while
$\Pinflamod$ will be the full subcategory of $\Lamod$ consisting of all
modules of finite projective dimension. 

In essence, our conceptual and notational backdrop is that developed in
successive steps in
\cite{\GePo, \Rin, \DoFr, \BuRi}, but some modifications to the
presentation of the relevant data will be more convenient for our purposes.
As is common, our notion of a `word' is based on the fixed presentation
$\la = K\Gamma / I$ as follows:  {\it Syllables\/} are elements of the set
$\Cal P \sqcup {\Cal P}^{-1}$, where
$\Cal P$ is the set of all paths in
$K\Gamma \setminus I$ and
${\Cal P}^{-1}$ consists of the formal inverses of the elements of $\Cal
P$.  The paths of length $0$, i.e., the vertices of $\Gamma$, are included
in $\Cal P$ and will be called the {\it trivial paths}; both these trivial
paths and their inverses are called {\it trivial syllables}.  For $p \in
\Cal P$, let
$(p^{-1})^{-1} = p$, so that
$({\Cal P}^{-1})^{-1} = \Cal P$.  {\it (Generalized) words\/} are
$\ZZ$-indexed sequences of pairs of syllables $w = (p_i^{-1} q_i)_{i \in
\ZZ}$ with $p_i, q_i \in
\Cal P$, which we also communicate as juxtapositions
$$\dots p_r^{-1} q_r \dots p_{-1}^{-1} q_{-1} p_0^{-1} q_0 p_1^{-1} q_1 \dots
p_s^{-1} q_s \dots$$ 
(note that syllables from $\Cal P$
alternate with syllables from
${\Cal P}^{-1}$) subject to the following constraints:
\roster
\item"$\bullet$"  For each $i \in \ZZ$, the starting points of $p_i$ and
$q_i$ coincide, but the first arrows of $p_i$ and $q_i$ are distinct
whenever both
$p_i$ and $q_i$ are nontrivial.
\item"$\bullet$"  For each $i \in \ZZ$, the end points of $q_i$ and
$p_{i+1}$ coincide, but the last arrows of $q_i$ and $p_{i+1}$ are distinct
whenever both $q_i$ and $p_{i+1}$ are nontrivial.
\item"$\bullet$" No trivial syllable occurs between two nontrivial syllables
(i\.e., the nontrivial syllables form a `connected component').
\endroster

A word  $w = (p_i^{-1} q_i)_{i \in \ZZ}$ will be called {\it finite}
in case, for all $i \gg 0$ and all $i \ll 0$, the syllables $p_i^{-1}$ (and
hence also the syllables
$q_i$) are trivial; finite words are also communicated as finite
juxtapositions
$(p_i^{-1} q_i)$ in which all nontrivial syllables are preserved.  More
generally, we do not insist on recording trivial syllables; keep in mind
that they can only occur at the left or right tail ends of a word.  It is
self-explanatory what we mean by a {\it left\/} or {\it right finite\/}
word.  Given a word $w = (p_i^{-1} q_i)_{i \in
\ZZ}$, the inverse of $w$ is defined as $w^{-1} = (q_{-i}^{-1} p_{-i})_{i \in
\ZZ}$ carrying the pair of syllables $q_{-i}^{-1} p_{-i}$ in position $i$.
With each word, we associate a (not necessarily finite) directed graph which
records the nontrivial syllables.  Namely, if $w = (p_i^{-1} q_i)_{i
\in
\ZZ}$ with
$p_i = a_{is(i)}
\dots a_{i1}$ and $q_i = b_{it(i)} \dots  b_{i1}$, where the $a_{ij}$ and
$b_{ij}$ are arrows, the graph of $w$ is

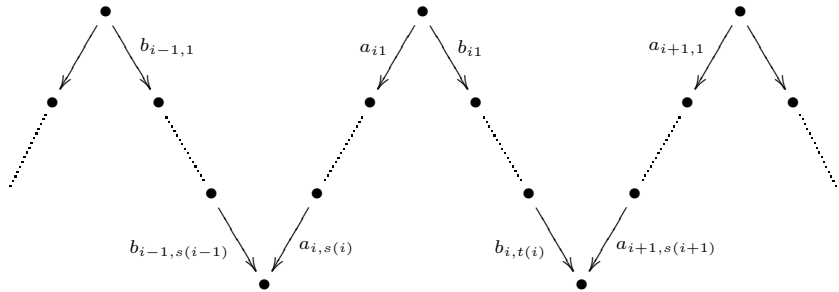
\ignore
$$\xymatrixrowsep{2pc}\xymatrixcolsep{0.75pc}
\xymatrix{
 &&\bullet \ar[dl] \ar[dr]^{b_{i-1,1}} &&&&&&\bullet \ar[dl]_{a_{i1}}
\ar[dr]^{b_{i1}} &&&&&&\bullet \ar[dl]_{a_{i+1,1}} \ar[dr] \\
 &\bullet \ar@{.}[dl] &&\bullet \ar@{.}[dr] &&&&\bullet \ar@{.}[dl]
&&\bullet \ar@{.}[dr] &&&&\bullet \ar@{.}[dl] &&\bullet \ar@{.}[dr] \\
 &&&&\bullet \ar[dr]_{b_{i-1,s(i-1)}} &&\bullet \ar[dl]^{a_{i,s(i)}}
&&&&\bullet \ar[dr]_{b_{i,t(i)}} &&\bullet \ar[dl]^{a_{i+1,s(i+1)}} &&&&\\
 &&&&&\bullet &&&&&&\bullet }$$
\endignore

\noindent  where the trivial syllables make no appearance, and the nodes are
identified with the primitive idempotents occurring as the starting and end
points of the arrows $a_{ij}$ and $b_{ij}$.  When less detail is required, a
simplified rendering of this graph will be preferred, namely:

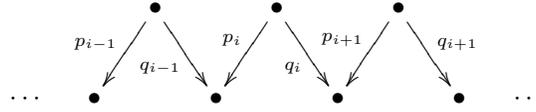
\ignore
$$\xymatrixrowsep{2pc}\xymatrixcolsep{1pc}
\xymatrix{
 &&\bullet \ar[dl]_{p_{i-1}} \ar[dr]_(0.55){q_{i-1}} &&\bullet
\ar[dl]_(0.45){p_i} \ar[dr]_(0.55){q_i} &&\bullet
\ar[dl]_(0.45){p_{i+1}} \ar[dr]^{q_{i+1}} \\
\cdots &\bullet &&\bullet &&\bullet &&\bullet &\cdots }$$
\endignore

\noindent  Graphs of words will only play a role in the proof of our main
theorem, while graphs of string and pseudo-band modules, as given below, will
be essential throughout our discussion.

To prepare for the upcoming definition, note that, for any nontrivial
two-syllable word $w = p^{-1}q$, there exists at most one arrow
$\alpha$ such that $(\alpha p)^{-1} q$ is again a word; analogously, there
exists at most one arrow $\beta$ making
$p^{-1}(\beta q)$ a word.

Each (generalized) word $w = (p_i^{-1} q_i)_{i \in \ZZ}$ gives rise to a
{\it (generalized) string module\/} $\St(w)$ defined as follows:   If $w$ is
trivial, say $w = e$, then $\St(w)$ is the simple module $\la e/ Je$.  Now
suppose that $w$ is nontrivial, and let $\supp(w)$ be the set of all those
integers $j$ for which either $p_j$ or $q_j$ is nontrivial. If, moreover,
the joint starting vertex of $p_i$ and $q_i$ is denoted by $e(i)$, then
$$\St(w) = \biggl(\, \bigoplus_{i \in \supp(w)}\la e(i) \biggr)
\biggm/ C \ , \quad \quad \text{where}$$
$$C\ \  = \ \biggl(\, \sum_{i, i+1 \in
\supp(w)} \la \bigl( q_i e(i) - p_{i+1} e(i+1) \bigr)\,\biggr)\ +\
C_{\text{left}} \ + \ C_{\text{right}},$$ with cyclic correction terms
$C_{\text{left}}$ and $C_{\text{right}}$ defined as follows:
$C_{\text{left}} = 0$ if either $\supp(w)$ is unbounded on the negative
$\ZZ$-axis or $l = \inf \supp(w)$ is an integer and
$(\alpha p_l)^{-1} q_l$ fails to be a word for all arrows
$\alpha$; in the remaining case, where $l \in \ZZ$ and there exists a
(necessarily unique) arrow $\alpha$ with the property that $(\alpha
p_l)^{-1} q_l$ is again a word, we set
$C_{\text{left}} = \la \alpha p_l e(l)$.  The right-hand correction term
$C_{\text{right}}$ is defined symmetrically.  Note that $\St(w)$ is finite
dimensional over $K$ precisely when $w$ is a finite word. Moreover, $\St(w)
\cong \St(w^{-1})$, a fact which allows us to pass back and forth between
$w$ and $w^{-1}$ as convenience dictates.  It is well-known that string
modules are indecomposable; in the finite dimensional case, this is proved in
\cite{\GePo, \Rin, \DoFr, \BuRi}, for infinite dimensional string modules in
\cite{\Kra}.

For $i \in \supp(w)$, let $x_i$ be the residue class of $e(i)$ in $\St(w)$
in the above presentation. Clearly, the family $(x_i)_{i\in \supp(w)}$
consists of top elements which generate $\St(w)$ and are
$K$-linearly independent modulo $J\St(w)$; by construction, they have the
property that $q_ix_i= p_{i+1}x_{i+1}$, whenever $q_i$ and $p_{i+1}$ are
both nontrivial. Any sequence of top elements of
$\St(w)$ with the listed properties is called a {\it standardized sequence
of top elements\/}. In the sense of \cite{\pre} and
\cite{\dom}, the module $\St(w)$ has a layered graph relative to any
standardized sequence of top elements:  It is the undirected variant of the
(directed) graph of $w$, layered in such a fashion that the vertices in the
$i$-th row from the top correspond to the simple composition factors in
$J^{i-1} \St(w)/ J^i \St(w)$.  We will usually indicate the chosen sequence
of top elements above the corresponding vertices in the first row of the
graph as illustrated below.

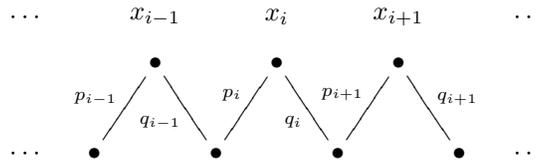
\ignore
$$\xymatrixrowsep{2pc}\xymatrixcolsep{1pc}
\xymatrix{
\dropup{4}{\cdots} &&\bullet \dropup{4}{x_{i-1}} \edge[dl]_{p_{i-1}}
\edge[dr]_(0.55){q_{i-1}} &&\bullet \dropup{4}{x_i}
\edge[dl]_(0.45){p_i} \edge[dr]_(0.55){q_i} &&\bullet \dropup{4}{x_{i+1}}
\edge[dl]_(0.45){p_{i+1}} \edge[dr]^{q_{i+1}} &&\dropup{4}{\cdots} \\
\cdots &\bullet &&\bullet &&\bullet &&\bullet &\cdots }$$
\endignore

The second class of representations of $\la$ which will be pivotal in our
discussion slightly generalizes the classical `band modules' (a
generalization which will prove convenient in the proof of the main
theorem).  This class consists entirely of finite dimensional
representations, but this time, they need not be indecomposable.  For a
description following the classical road, suppose that
$v= p_0^{-1}q_0
\dots p_t^{-1}q_t$ is a finite word with $t\ge 0$ and $p_0$, $q_t$ both
nontrivial; by our conventions, this amounts to the same as to require that
all $p_i$ and $q_i$ be nontrivial. We call
$v$ {\it primitive\/} if
\roster
\item"$\bullet$" the juxtaposition $v^2 =vv$ is again a word (in which case
all powers
$v^r$ are words), and
\item"$\bullet$" $v$ is not itself a power of a strictly shorter word.
\endroster

In addition to the primitive word $v$,  let $r$ be a positive integer and
$\phi: K^r \rightarrow K^r$ a cyclic automorphism (meaning that $\phi$ turns
$K^r$ into a cyclic $K[X]$-module)  with Frobenius companion matrix
$$\pmatrix 0&\cdots&0&c_1\\ 1&\ddots &&\vdots\\ &\ddots&0&\vdots\\
0&\cdots&1&c_r \endpmatrix.$$
 Then the {\it pseudo-band module\/}
$\Bd(v^r,\phi)$ is defined as follows: Let $x_{10},\dots,x_{1t}$,
$x_{20}, \dots, x_{2t}$, $\dots, x_{r0}, \dots, x_{rt}$ be the standardized
sequence of top elements of $\St(v^r)$ following the definition of a string
module; in particular, $q_jx_{ij}= p_{j+1}x_{i,j+1}$ for $0\le i\le r$ and
$1\le j<t$, and
$q_tx_{it}= p_0x_{i+1,0}$ for $i<r$. Then
$$\Bd(v^r,\phi)= \St(v^r)/ \la \bigl( q_rx_{rt}- \sum_{i=1}^r c_ip_0x_{i0}
\bigr)$$ by definition. If the residue classes of the $x_{ij}$ in
$\Bd(v^r,\phi)$ are denoted by $y_{ij}$, then the latter are top elements of
$\Bd(v^r,\phi)$ which are
$K$-linearly independent modulo the radical, generate $\Bd(v^r,\phi)$, and
satisfy the above equations, as well as the additional one
$$q_ty_{rt}= \sum_{i=1}^r c_i p_0y_{i0}.$$ Any sequence of top elements with
these properties is, in turn, called a standardized sequence of top elements
of the pseudo-band module $\Bd(v^r,\phi)$. Relative to such a sequence,
$\Bd(v^r,\phi)$ can be depicted in the form

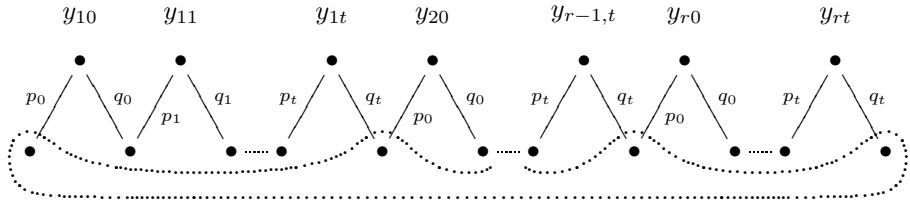
\ignore
$$\xymatrixrowsep{2pc}\xymatrixcolsep{0.67pc}
\xymatrix{
 &\bullet \dropup{4}{y_{10}} \edge[dl]_{p_0} \edge[dr]^{q_0} &&\bullet
\dropup{4}{y_{11}} \edge[dl]^(0.55){p_1}
\edge[dr]^{q_1} &&&\bullet \dropup{4}{y_{1t}}
\edge[dl]_{p_t} \edge[dr]^{q_t} &&\bullet \dropup{4}{y_{20}}
\edge[dl]^(0.55){p_0}
\edge[dr]^{q_0} &&&\bullet \dropup{4}{y_{r-1,t}}
\edge[dl]_{p_t} \edge[dr]^{q_t} &&\bullet \dropup{4}{y_{r0}}
\edge[dl]^(0.55){p_0}
\edge[dr]^{q_0} &&&\bullet  \dropup{4}{y_{rt}}
\edge[dl]_{p_t} \edge[dr]^{q_t}\\
\bullet &&\bullet &&\bullet  \ar@{.}[r] &\bullet &&\bullet &&\bullet
\ar@{.}[r] &\bullet
\save[0,0]+(-1,-2);[0,-1]+(1,-2) **\crv{~*=<2.5pt>{.} [0,0]+(0,-3)
&[0,1]+(0,-3) &[0,2]+(-2,3) &[0,2]+(2,3) &[0,3]+(0,-3) &[0,6]+(0,-3)
&[0,7]+(-2,3) &[0,7]+(3,3) &[0,7]+(3,-6) &[0,6]+(0,-6) &[0,-9]+(0,-6)
&[0,-10]+(-3,-6) &[0,-10]+(-3,3) &[0,-10]+(2,3) &[0,-9]+(0,-3)
&[0,-4]+(0,-3) &[0,-3]+(-2,3) &[0,-3]+(2,3) &[0,-2]+(0,-3) &[0,-1]+(0,-3)
}\restore
 &&\bullet &&\bullet  \ar@{.}[r] &\bullet &&\bullet }$$
\endignore
\vskip0.3truein

\noindent where the dotted line that encircles the vertices standing for the
elements $q_ty_{rt}$, $p_0 y_{i0},\dots,
\allowmathbreak p_0y_{r0}$ in the socle of $\Bd(v^r,\phi)$ indicates that
the $K$-space spanned by these $r+1$ elements has dimension $r$ only.

In case the automorphism $\phi$ of $K^r$ is irreducible, the pseudo-band
module $\Bd(v^r,\phi)$ is called a {\it band module\/}.  It is readily seen
that a pseudo-band module is a band module if and only if it is
indecomposable (see, e\.g.,
\cite{\BuRi}). Note that, in contrast to the string case, the graph of a
band module $\Bd(v^r,\phi)$ does not pin down the latter up to isomorphism,
unless the scalars $c_1,\dots,c_r$ are recorded. Moreover, observe that, if
we subject the pairs of syllables $p_i^{-1}q_i$ of the underlying primitive
word $v$ to a cyclical permutation resulting in $\vhat$ say, then
$\Bd(v^r,\phi)
\cong \Bd(\vhat^r, \phi)$.

The pivotal role played by the finite dimensional string and band modules is
apparent from the following classification result, which will be used
extensively in the sequel. In its present form, it was established by Butler
and Ringel, but the ideas go back to Gelfand and Ponomarev who determined
the finite dimensional representation theory of a somewhat more restricted
class of algebras.

\proclaim{Theorem 0} {\rm (See \cite{\GePo, \Rin, \Bon, \DoFr,
\BuRi})} The finitely generated string and band modules are precisely the
indecomposable objects of $\lamod$. \qed\endproclaim
\medskip

\head 3.  Syzygies and the homological dimensions of string algebras
\endhead

This short section provides a first installment of evidence that, not only
from a represen\-ta\-tion-the\-or\-et\-ic, but also from a homological
viewpoint, string algebras show `tame behavior'.  Not only can the global
dimension of a string algebra be computed algorithmically from quiver and
relations as we will shortly see, but this is true more generally for its
finitistic dimensions.

Our first proposition determines the syzygies of string and band modules.
All of these are direct sums of cyclic string modules which can be
described in terms of the string and band data.  Since we will repeatedly
invest this information in subsequent sections, we will be explicit,
starting with the slightly cumbersome notation required to pin down syzygies.
Suppose that $p$ and
$q$ are paths, not both trivial, such that
$p^{-1}q$ is a word, i.e., $p$ and $q$ both start in the same vertex $e$,
and if both of these paths are nontrivial, they have distinct first arrows. 
Then, clearly, the string module
$\St(p^{-1}q)$ is a factor module of $\la e$.  Let $\fp$ and $\fq$ be the
unique paths starting in the end points of $p$ and $q$, respectively, with
the property that $\St \bigl( (\fp p)^{-1}(\fq q) \bigr) \cong \la e$.  In
other words, if $p$ is nontrivial, then $\fp$ is the longest path such that
$\fp p \in K\Gamma
\setminus I$, the path $\fq$ having an analogous description if
$q$ is nontrivial;  if on the other hand, $p$ is trivial, then $q$ is not,
and $\fp$ is the longest path starting in $e$ which does not contain the
first arrow of $q$ as a right subpath (in particular, $\fp = e$ in case $q$
is nontrivial and
$\la e$ is uniserial).  Finally, given any nontrivial path $u$, let
$u^{(0)}$ be the path of length
$\ge 0$ obtained from $u$ through deletion of the first arrow, and set
$u^{(0)} = 0 \in K\Gamma$ if $u$ is trivial.

We observe that the first syzygy of any cyclic string module $\St(p^{-1}q)$
equals
$\St(\fp^{(0)}) \oplus \St(\fq^{(0)})$ if $p$, $q$ are not both trivial, and
equals
$Je$ if $p=q=e$.  Less obvious situations are addressed in the first
proposition, the proof of which is immediate from the definitions and the
graphical methods for determining syzygies developed in \cite{\pre}.

\proclaim{Proposition 1}

{\rm (1)(a)}  Suppose $w$ is a finite nontrivial word of the form
$p_0^{-1} q_0 \dots  p_t^{-1} q_t$ with $q_0$ and $p_t$ nontrivial.  Then
the first syzygy of the finite dimensional string module $\St(w)$ is

$$\St(\fp_0^{(0)}) \oplus \bigoplus_{i=0}^{t-1} \St(\fp_{i+1}^{-1}
\fq_i) \oplus \St(\fq_t^{(0)}),$$ where the paths $\fp_i$, $\fq_i$,
$\fq_0^{(0)}$, and $\fq_t^{(0)}$ are as introduced above.

{\rm (b)} Now suppose that $w = \dots p_{-1}^{-1} q_{-1} p_0^{-1} q_0
p_1^{-1} q_1 \dots$ is a word with $p_i$, $q_i$ nontrivial for all $i \in
\ZZ$.  The first syzygy of the infinite dimensional string module $\St(w)$
equals
$$\bigoplus_{i \in \ZZ} \St(\fp_{i+1}^{-1}\fq_i).$$
\medskip

{\rm (2)}  If $v = p_0^{-1} q_0 \dots  p_t^{-1} q_t$ is a primitive word
with $p_0$ and $q_t$ nontrivial, $r$ a positive integer, and $\phi$ a cyclic
automorphism of $K^r$, then the first syzygy of the pseudo-band module
$\Bd(v^r, \phi)$ is
$$\bigoplus_{i=0}^{t-1} \bigl(\St(\fp_{i+1}^{-1} \fq_i)
\bigr)^r \oplus \bigl( \St(\fp_0^{-1} \fq_t) \bigr)^r.$$

Furthermore, the following statements are equivalent:
\roster
\item"(i)" For {\rm some} positive integer $r$ and {\rm some} cyclic
automorphism $\phi$ of $K^r$, the pseudo-band module $\Bd(v^r, \phi)$
belongs to $\pinflamod$.
\item"(ii)" For {\rm any} positive integer $r$ and {\rm any} cyclic
automorphism $\phi$ of $K^r$, the pseudo-band module $\Bd(v^r, \phi)$
belongs to $\pinflamod$.
\item"(iii)" The generalized string module $\St(\dots vvv\dots)$ belongs to
$\Pinflamod$
 \qed \endroster
\endproclaim

In view of Theorem 0 (Section 2), we glean from Proposition 1 that the
syzygies of any finitely generated
$\la$-module $M$ are direct sums of cyclic string modules.  Next we will see
that this remains true even when we drop the requirement that $M$ be
finitely generated.

\proclaim{Proposition 2}  Every submodule of a projective left
$\la$-module is a direct sum of string modules of the form
$\St(p^{-1}q)$, where $p$ and $q$ are paths.

In particular, syzygies of arbitrary $\la$-modules are direct sums of cyclic
string modules which embed in $J$, and all second syzygies are direct sums
of uniserial modules.
\endproclaim

\demo{Proof}  To verify the first assertion, let $P$ be a projective left
$\la$-module and
$M
\subseteq P$ a submodule.  Since all of the indecomposable projective left
$\la$-modules are cyclic string modules, it is harmless to assume that $M
\subseteq JP$.

Write $M$ as a directed union of finitely generated submodules, say
$M = \bigcup_{i \in I} M_i$.  Then all $M_i$ are syzygies of finitely
generated modules and, by the preceding remark, the
$M_i$ are direct sums of string modules with simple tops.  Since (up to
isomorphism) there are only finitely many string modules of of the form
$\St(p^{-1}q)$, this entails that the direct sum
$\bigoplus_{i \in I} M_i$ is $\Sigma$-algebraically compact.  Therefore, by
 \cite{\HuSm, Observation 3.1} or \cite{\Len}, the category $\Add
\bigl(\bigoplus_{i \in I} M_i \bigr)$ of arbitrary direct sums of direct summands
of the
$M_i$ is closed under direct limits.  In particular,
$M$ is in turn a direct sum of cyclic string modules, as claimed.

The final statements are immediate consequences.\qed \enddemo

Recall that the left little and big finitistic dimensions of any finite
dimensional algebra $\Delta$ are defined as
$$\align \lfindim \Delta &= \sup\{\pdim M \mid M \in \pinfdeltamod \} \\
\lFindim \Delta &= \sup\{\pdim M \mid M \in \Pinfdeltamod \}, \endalign$$
respectively.  According to \cite{\pre} or \cite{\dom}, the finitistic
dimensions of any monomial relation algebra can be computed up to an error of
$1$ by means of a simple graphical method.  However, even for monomial
relation algebras, the little finitistic dimension may be strictly smaller
than the big finitistic dimension \cite{\dom}, and, for more general finite
dimensional algebras, the difference $\lFindim \Delta -
\lfindim \Delta$ attains arbitrarily high values in $\NN$ \cite{\Sma}.  Our
first theorem excludes such `pathologies' for string algebras and pins down
the finitistic dimensions in terms of the cyclic string modules of finite
projective dimension; the latter are finite in number and easy to construct from
$\Gamma$ and $I$.  To that end, consider the  set
$\Cal T$ of those cyclic string modules in
$\pinflamod$ which are contained in the radical $J$ of $\la$, that is,
$$\Cal T = \{\St(p^{-1}q) \in \pinflamod \mid \St(p^{-1}q) \text{\ embeds
in\ } J \};$$ here the paths $p$ and $q$ may be trivial, and so, in
particular, $\Cal T$ includes all uniserial left modules of finite
projective dimension contained in $J$.  The proof of Theorem 3 is an
immediate consequence of Proposition 2.

\proclaim{Theorem 3} $\lfindim \la = \lFindim \la = t+1$, where $t =
\sup\{\pdim M \mid M \in \Cal T\}$ in case $\Cal T$ is nonempty, and $t=-1$
otherwise. \qed
\endproclaim

The following example shows that shrinking $\Cal T$ to the set of all {\it
uniserial\/} left modules from $\pinflamod$ contained in $J$ does not leave
the conclusion of Theorem 3 intact.

\example{Example 4}  Let $\la = K\Gamma / I$, where $\Gamma$ is the quiver

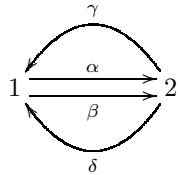
\ignore
$$\xymatrixcolsep{4pc}
\xymatrix{ 1 \ar[r]<0.75ex>^{\alpha} \ar[r]<-0.75ex>_{\beta} &2
\ar@/_2pc/[l]_{\gamma} \ar@/^2pc/[l]^{\delta}
 }$$
\endignore

\noindent and $I \subseteq K\Gamma$ is the monomial ideal with the property
that the indecomposable projective left
$\la$-modules have graphs

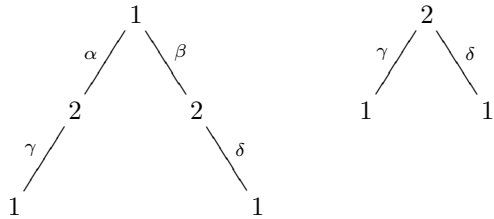
\ignore
$$\xymatrixrowsep{2pc}\xymatrixcolsep{1pc}
\xymatrix{
 &&1 \edge[dl]_{\alpha} \edge[dr]^{\beta} &&&&&2 \edge[dl]_{\gamma}
\edge[dr]^{\delta}\\
 &2 \edge[dl]_{\gamma} &&2 \edge[dr]^{\delta} &&&1 &&1\\ 
1 &&&&1 }$$
\endignore

Then $\Cal T$ is the singleton containing
$\la e_2 \cong \la(\alpha - \beta)$, whence, by Theorem 3, the left little
and big finitistic dimensions of $\la$ are equal to
$1$. 
\endexample

\head 4. Background on contravariant finiteness and phantoms \endhead

We give a brief summary of that part of the theory of contravariant
finiteness of a full subcategory $\A \subseteq
\lamod$ which will be relevant here.  The two categories on which we will
focus in the sequel are $\pinflamod$ and $\sinflamod$, the full subcategory
of $\pinflamod$ consisting of the finite direct sums of string modules; note
that the latter category depends, if only to a `minor' degree, on the
coordinatization $(\Gamma, I)$ of the string algebra
$\la$. 

Since the background material we need is quite general, we let
$\Delta$ be any finite dimensional algebra for the moment, and
$\A \subseteq \deltamod$ a full subcategory which is closed under finite
direct sums and direct summands.  Recall from  \cite{\AuSm} that $\A$ is said
to be {\it contravariantly finite} in
$\deltamod$ in case each module $M \in \deltamod$ has a (right) $\A$-{\it
approximation}, that is a homomorphism
$\varphi: A
\rightarrow M$ with $A \in \A$ such that the induced sequence of functors
$$\Hom_{\Delta}(-,A)|_{\A} \longrightarrow \Hom_{\Delta}(-,M)|_{\A}
\longrightarrow 0$$  is exact; in other words, the latter says that every
map in $\Hom_{\la}(B,M)$ with $B \in \A$ factors through $\varphi$. In the
sequel, we will suppress the qualifier `right'.  Provided that $M$ has an
$\A$-approximation, there is a {\it minimal\/} such approximation which
embeds, as a direct summand, into all other $\A$-approximations of $M$.  In
particular, such a  minimal $\A$-approximation of $M$ is unique up to
isomorphism; it is therefore only a mild abuse of language to refer to it as
{\it the minimal $\A$-approximation of $M$}.  In case $\A = \pinfdeltamod$,
the existence of approximations for the simple left $\Delta$-modules $S_i$,
$1 \le i \le n$, already guarantees contravariant finiteness of
$\pinfdeltamod$.  Moreover, in case of existence, the minimal
$\pinflamod$-approximations
$A_i$ of the $S_i$ impinge on the structure of an arbitrary object in
$\pinflamod$ as follows (see \cite{\AuRe}):  A finitely generated
$\Delta$-module $M$ has finite projective dimension precisely when it is a
direct summand of a module $X$ that has a finite filtration with successive
factors in $\{A_1, \dots, A_n\}$.  This structure theory was extended to
$\Pinfdeltamod$, the category of {\it all\/} left $\Delta$-modules of finite
projective dimension by the authors in \cite{\HuSm}:  If
$\pinfdeltamod$ is contravariantly finite, the objects of the `big' category
$\Pinfdeltamod$ are precisely those which are direct limits of objects $X$
having filtrations $X = X_0 \supseteq X_1 \supseteq \cdots \supseteq X_m =
0$ with $X_i/X_{i+1} \cong A_{k(i)}$ for all $i < m$; in particular, the big
and little left finitistic dimensions of $\Delta$ coincide in this situation
and are attained on
$\{A_1, \dots, A_n\}$. This illustrates the pivotal role played by the
minimal approximations of the simple modules in case
$\pinflamod$ is contravariantly finite.

$\A$-{\it phantoms\/} were introduced in \cite{\HaHZ},
originally mainly with the aim of providing criteria for contravariant
finiteness of
$\A$.  Here they serve a dual purpose:  the one just mentioned and that
of replacing minimal $\A$-approximations in case such approximations fail to
exist.  Our phantoms under (3) below are not to be confused with the phantom
maps of algebraic topology, or the topology-inspired phantom maps in the
modular representation theory of groups, as introduced by Benson and
Gnacadja in
\cite{\BenGna}.

\definition{Definition} Assume $\A \subseteq \deltamod$ to be closed under
direct summands and finite direct sums.  Moreover, fix
$M \in \deltamod$.

{\rm (1)}  Let $\C$ be a subcategory of $\A$.  A {\it relative
$\C$-approximation of $M$ in\/} $\A$ is a homomorphism $f: A
\rightarrow M$ with $A \in \A$ such that all maps in
$\Hom_{\Delta}(\C,M)$ factor through $f$.  In this situation, we will also
refer to the module $A$ as a relative $\C$-approximation of $M$ in $\A$.

If $\C = \{C\}$, we write `relative $C$-approximation' for short.

{\rm (2)}  A finitely generated module $H$ is an $\A$-{\it phantom of\/} $M$
in case there exists an object $C$ in $\A$ with the property that $H$ occurs
as a subfactor of {\it every} relative
$C$-approximation of $M$ in $\A$.

More generally, a module $H \in \Deltamod$ will be called an
$\A$-{\it phantom\/} of $M$ if each of its finitely generated submodules is
a phantom in the sense just defined.

{\rm (3)}  Given an $\A$-phantom $H$ of $M$ and a subcategory $\C$ of $\A$,
a homomorphism $\varphi: H \rightarrow M$ is called an {\it effective
$\C$-phantom\/} of $M$, provided that $H$ is a direct limit of objects in
$\C$, and every map in
$\Hom_{\Delta}(\C,M)$ factors through $\varphi$.  In that case, we will also
say that $H$ is an effective $\C$-phantom of $M$.
\enddefinition

We add a few comments to set up an intuitive backdrop for the concept of a
phantom.  The terminology is to evoke a `phantom image' of an object,
assembled from witness reports, to aid a search effort. In that spirit,
a finitely generated module $H$ is an
$\A$-phantom of $M$ precisely when there is a witness $C \in \A$ testifying
to the effect that the source of any homomorphism in
$\Hom_{\la}(\A, M)$, which permits factorization of all maps in
$\Hom_{\la}(C, M)$, has an epimorphic image containing $H$.

In case $\A$-approximations of $M$ exist, the minimal one, $A$ say, is the
only effective $\A$-phantom of $M$, and the class of all $\A$-phantoms of $M$
coincides with the class of subfactors of
$A$. So, in this situation, constructing $\A$-phantoms of
$M$ amounts to assembling `phantom images' of the module $A$ that one would
like to track down;  existence is recognized in the process if
one can argue that the $K$-dimensions of such phantoms need to be bounded (cf\.
the existence theorem below).  Otherwise, namely when
$M$ fails to have
$\A$-approximations, $\A$-phantoms of $M$ still provide minimal building
blocks of objects through which all maps in
$\Hom_{\la}(\A, M)$ can be factored. Of course, effective
phantoms of $M$ hold the highest structural interest also in this case: Effective
$\C$-phantoms, where $\C \subseteq \A$, take over the role of minimal
approximations, in that they carry full complements of information on how
$\C$ relates to $M$.

Note that, for any finite subcategory $\C$ of $\A$, relative
$\C$-approximations of $M$ in $\A$ exist; to obtain candidates, we only have
to add up a sufficient number of copies of the objects in $\C$.  Moreover,
observe that, given a subclass $\D \subseteq
\C$, every relative $\C$-approximation of $M$ in $\A$ is also a relative
$\D$-approximation. Hence calling for an object $C$ in
$\A$ such that $H$ is a subfactor of {\it every\/} relative
$C$-approximation of $M$ in $\A$  -- as we do in the definition of a
finitely generated phantom  --  places strong `minimality pressure' on $H$. 
On the other hand, the class of $\A$-phantoms of $M$ is closed under
subfactors and direct limits of directed systems, which often makes it
enormous.  This slack in the definition facilitates the search for phantoms
and thus makes them an expedient tool in proving failure of contravariant
finiteness: Indeed, the existence of $\A$-phantoms of unbounded lengths of a
given module $M$ signals non-existence of an $\A$-approximation of $M$. We
conclude this sketch with the following existence result.

\proclaim{Theorem} {\rm (see \cite{\HaHZ})} For $M$ in $\deltamod$ and $\A
\subseteq \deltamod$ a subcategory as above, the following conditions are
equivalent:

{\rm (1)}  $M$ fails to have an $\A$-approximation.

{\rm (2)}  $M$ has $\A$-phantoms of arbitrarily high finite $K$-dimensions.

{\rm (3)}  $M$ has an $\A$-phantom of infinite $K$-dimension.

{\rm (4)}  There exists a countable subclass $\C \subseteq \A$ such that
$M$ has an effective $\C$-phantom of infinite $K$-dimension. \qed
\endproclaim

 For a more extensive overview of contravariant finiteness results, we refer
the reader to
\cite{\Huiathens}.

\head 5.  Statement of the main result, consequences, and examples \endhead

Given a generalized word $w = (p_i^{-1} q_i)_{i \in \ZZ}$, which is either
trivial or has the property that $\length(p_0) +
\length(q_0) > 0$, we will refer to the joint starting point of
$p_0$ and $q_0$ as the {\it center} of $w$; denote this center by
$e$. If we wish to emphasize the centered viewpoint, we will also refer to
$w$ as a {\it generalized word centered at\/} $e$.  Each generalized word
$w$ centered at $e$ comes paired with an obvious homomorphism
$\varphi:\St(w) \rightarrow \la e / Je$:  Indeed, consider the standard
presentation
$$\St(w) = \biggl( \bigoplus_{i \in \supp(w)} \la e(i) \biggr) \biggm/\ C$$
as specified in Section 2; here $e(i)$ is again the starting vertex of $p_i$,
$q_i$ for $i \in \ZZ$.  Let $x_i \in \St(w)$ be the residue class of $e(i)$.
Then there exists a unique homomorphism $\varphi: \St(w) \rightarrow \la
e/Je$ with
$\varphi(x_0) = e + Je$ and $\varphi(x_i) = 0$ for $|i| \ge 1$; we will
refer to it as the {\it canonical map} of the centered word
$w$.

Moreover, we call a generalized word $w$ {\it left periodic}, resp\. {\it
right periodic}, in case $w = \dots uuu w_1$, resp\.,
$w = w_2 vvv \dots$, with $u$, $v$ either trivial or primitive, and $w_1$,
$w_2$ left, resp\. right, finite. (On the side, we point out that the set of
those left and right periodic words which are twosided infinite coincides
with the union of the `periodic' and `biperiodic $\ZZ$-words' introduced by
Ringel in \cite{\Rintwo}, while our left finite and right periodic words are
`periodic' or `almost periodic $\NN$-words' in Ringel's terminology.)

Finally, we recall that $\sinflamod$ denotes the full subcategory of
$\pinflamod$ having as objects all finite direct sums of string modules of
finite projective dimension.

We are now in a position to state our main theorem.

\proclaim{Theorem 5}  As before, let $\la = K\Gamma / I$ be a string algebra
with simple left modules $S_i = \la e_i / J e_i$,
$1 \le i \le n$.  Then there exist centered generalized words $w_i =
w(S_i)$, unique up to inversion, with  the following properties:
\smallskip

{\rm (I)}  Each $w_i$ is centered at $e_i$, left and right periodic, and  can
be effectively constructed from $\Gamma$ and $I$ in a number of steps which
grows polynomially with $\dim_K \la$.
\smallskip

{\rm (II)}  Each of the generalized string modules $\St(w_i)$ has finite
projective dimension, and the canonical map
$\varphi_i:
\St(w_i)
\rightarrow S_i$ is an {\rm effective}
$\sinflamod$-phantom of $S_i$.
\smallskip

{\rm (III)}  The category $\pinflamod$ is contravariantly finite in $\lamod$
if and only if the words $w_1, \dots, w_n$ are all finite.  In the positive
case, the canonical maps $\varphi_i:
\St(w_i) \rightarrow S_i$ are the minimal
$\pinflamod$-approximations of the simple modules.  More precisely, for each
$i \in \{1, \dots, n\}$, the following conditions are equivalent:

\roster
\item"(i)" $w_i$ is finite.

\item"(ii)"  $S_i$ has a $\pinflamod$-approximation.

\item"(iii)"  $\varphi_i: \St(w_i)
\rightarrow S_i$ is the minimal $\pinflamod$-approximation of $S_i$.

\item"(iv)" $S_i$ has an $\sinflamod$-approximation.

\item"(v)" $\varphi_i: \St(w_i)
\rightarrow S_i$ is the minimal $\sinflamod$-approximation of $S_i$.
\endroster

\noindent  Finally, if the equivalent conditions {\rm (i) -- (v)} are
satisfied, the top of $\St(w_i)$ has dimension at most $4n$.
\endproclaim

A proof will be given in Section 7.  As will be further substantiated when
we describe and explore the generalized words
$w_i$ of Theorem 5, they encode essentially all of the information required
to understand the homology of $\la$.  Without providing an explicit
algorithm, we mention that, in case of contravariant finiteness of
$\pinflamod$ in $\lamod$, the minimal
$\pinflamod$-approximations of arbitrary string and band modules can readily
be constructed from the $\St(w_i)$.

Theorem 5  remains unaffected if we replace the words $w_i$ by their
inverses; indeed, $\St(w) \cong \St(w^{-1})$ for any generalized word $w$,
and if $w$ is centered at $e$, the obvious isomorphism, namely the flip
about the center, preserves the center  and takes the canonical map $\St(w)
\rightarrow \la e/ Je$ to the canonical map $\St(w^{-1}) \rightarrow \la e/
Je$. Consequently, we will invert whenever convenient.  Note, moreover, that
Theorem 5 guarantees uniqueness of any of the string modules
$\St(w_i)$ up to isomorphism, whenever $w_i$ is finite.  On the other hand,
the isomorphism type of $\St(w_i)$ may depend on the coordinatization in
general, a fact which reflects the dependence of $\sinflamod$ on the
coordinatization;  see Section 8, Example 23.   Yet, this failure of
uniqueness is `minor', in that the graph of $\St(w_i)$, minus the labeling
of the edges, is invariant up to a flip about the central axis; this is a
consequence of Proposition 16 below.  In categorical terms, if $w_1, \dots, w_n$
and $w_1', \dots, w_n'$ are centered words having the properties described in
Theorem 5 relative to two eligible coordinatizations of our string algebra
$\la$, there exists a Morita self-equivalence $F: \Lamod
\rightarrow \Lamod$ such that $\St(w_i') \cong F \bigl(
\St(w_i)\bigr)$ for all $i$, and $F$ carries the canonical epimorphisms
$\varphi_i: \St(w_i) \rightarrow S_i$ to those of the centered words $w_i'$.

 We conclude the section with two immediate consequences of the main
theorem, followed by examples showing that neither can be extended to
arbitrary special biserial algebras.

\proclaim{Corollary 6}  Suppose that $\la$ is a string algebra. If
$\pinflamod$ is contravariantly finite, the minimal
$\pinflamod$-ap\-prox\-i\-ma\-tions of the simple left
$\la$-modules are string modules.  In particular, they are indecomposable.
\qed \endproclaim

While, under the hypothesis of the corollary, the category  of all finite
direct sums of string modules is always closed with respect to minimal
$\pinflamod$-approximations (we do not include a proof for this fact),
indecomposability is not preserved in passing from a string module to its
minimal approximation in general.

Actually, the property of being a string module may appear somewhat
artificial, since it usually depends on the
coordinatization of
$\la$.  However, in light of Theorems 3 and 5, it becomes apparent that the
concepts of `string' and `band' are more than devices permitting an
explicit classification of the finitely generated indecomposable
representations from quiver and relations of $\la$.  Indeed, we see that the
class of all string modules determines (irrespective of the chosen
coordinatization) the homological properties of the category
$\lamod$.

Since the base field $K$ does not enter into the structure of string modules
and their syzygies, this, in turn, guarantees that the homology of string 
algebras is a purely combinatorial game which is governed by the graphs of
the indecomposable projective modules alone.  In particular, the various
homological dimensions of a string algebra $\la = K\Gamma / I$ are
completely determined by the quiver $\Gamma$ and any set of paths generating
the ideal
$I$.  To see that this does not extend to arbitrary monomial relation
algebras, compare, e\.g\., \cite{\dep}.

\proclaim{Corollary 7}  The class of string modules determines the
homological dimensions of a string algebra $\la \allowmathbreak =
\allowmathbreak K\Gamma/I$, as well as the contravariant finiteness status
of the subcategory $\pinflamod$.  In particular, these data depend only on
the quiver $\Gamma$ and the paths in
$I$, not on the base field. \qed
\endproclaim

Both corollaries fail for special biserial algebras in general.

\example{Example 8} We present a finite dimensional special biserial algebra
$\Delta$ and a simple
$\Delta$-module $S_1$ whose $\pinfdeltamod$-approximation splits into three
nontrivial summands.  Suppose that
$\Delta$ is a path algebra modulo relations over $K$ with $\dim_K \Delta /
J(\Delta) = 8$ such that the indecomposable projective left modules have the
following graphs:

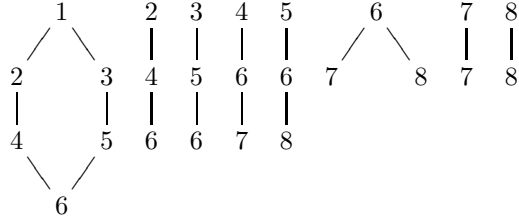
\ignore
$$\xymatrixrowsep{1pc}\xymatrixcolsep{0.5pc}
\xymatrix{
 &1 \edge[dl] \edge[dr] &&2 \edge[d] &3 \edge[d] &4 \edge[d] &5 \edge[d] &&6
\edge[dl]
\edge[dr] &&7 \edge[d] &8 \edge[d]  \\ 2 \edge[d] &&3 \edge[d] &4 \edge[d]
&5 \edge[d] &6 \edge[d] &6 \edge[d] &7 &&8 &7 &8
\\ 4 \edge[dr] &&5 \edge[dl] &6 &6 &7 &8 \\
 &6 }$$
\endignore

\noindent Clearly, $\Delta$ is a special biserial algebra. Moreover, it is
not difficult to check that the minimal
$\pinfdeltamod$-approxi\-ma\-tion of $S_1$ is the module with graph

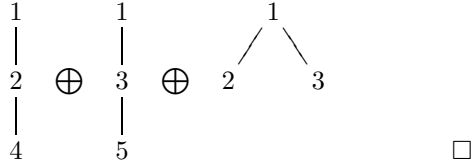
\ignore
$$\xymatrixrowsep{1pc}\xymatrixcolsep{0.5pc}
\xymatrix{ 1 \edge[d] &&1 \edge[d] &&&1 \edge[dl] \edge[dr] \\ 2 \edge[d]
&\bigoplus &3 \edge[d] &\bigoplus &2 &&3  \\ 4 &&5 &&&& &&&&\square }$$
\endignore
\endexample

\example{Example 9} For the following finite dimensional special biserial
algebra $\Delta$, the category $\sinfdeltamod$ has finite representation
type and is thus contravariantly finite  --  in particular, all simple
modules have $\sinflamod$-approximations  --  whereas $\pinfdeltamod$ fails
to be contravariantly finite in
$\deltamod$.  Let $\Delta$ be a path algebra with relations such that the
graphs of the indecomposable projective left
$\Delta$-modules are

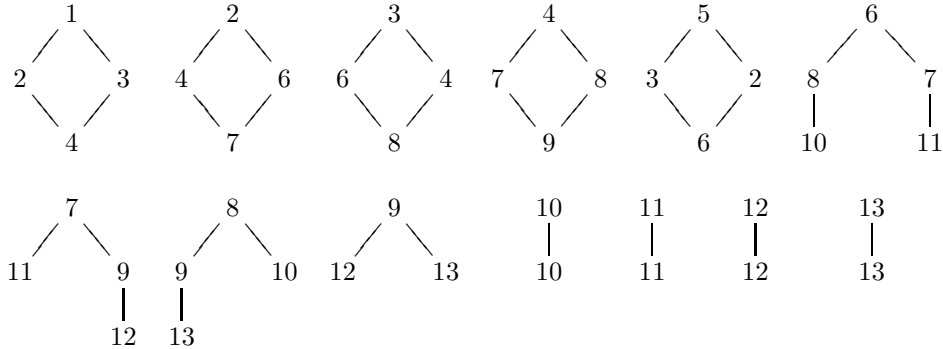
\ignore
$$\xymatrixrowsep{1pc}\xymatrixcolsep{0.5pc}
\xymatrix{
 &1 \edge[dl] \edge[dr] &&&2 \edge[dl] \edge[dr] &&&3 \edge[dl] \edge[dr]
&&&4 \edge[dl] \edge[dr] &&&5 \edge[dl] \edge[dr] &&&6 \edge[dl] \edge[dr]
\\ 2 \edge[dr] &&3 \edge[dl] &4 \edge[dr] &&6 \edge[dl] &6 \edge[dr] &&4
\edge[dl] &7 \edge[dr] &&8 \edge[dl] &3 \edge[dr] &&2 \edge[dl] &8 \edge[d]
&&7 \edge[d] \\
 &4 &&&7 &&&8 &&&9 &&&6 &&10 &&11 \\
 &7 \edge[dl] \edge[dr] &&&8 \edge[dl] \edge[dr] &&&9 \edge[dl] \edge[dr]
&&&10 \edge[d] &&11 \edge[d] &&12 \edge[d] &&13 \edge[d] \\ 11 &&9 \edge[d]
&9 \edge[d] &&10 &12 &&13 &&10 &&11 &&12 &&13 \\
 &&12 &13 }$$
\endignore

\noindent  It is straightforward (if somewhat tedious) to check that the
only finitely generated string modules of finite projective dimension in
$\deltamod$ are the indecomposable projective modules $\Delta e_6$ through
$\Delta e_{13}$, and the module with graph \ignore $\vcenter{
\xymatrixrowsep{0.6pc}\xymatrixcolsep{0.3pc} \xymatrix{ &7
\edge[dl] \edge[dr] &&8 \edge[dl] \edge[dr] \\ 11 &&9 &&10 } }$
\endignore.  This shows that $\sinfdeltamod$ has finite type. Note, in
particular, that there are no string modules of finite projective dimension
which have a copy of $S_1$ in their tops, whence the zero map $0 \rightarrow
S_1$ is the minimal
$\sinflamod$-approximation of $S_1$.  On the other hand, $S_1$ fails to have
a $\pinfdeltamod$-approxi\-ma\-tion;  indeed, since each band module
$\Bd(v^r,\phi)$, where $v$ is the primitive word with graph \ignore
$\vcenter{
\xymatrixrowsep{0.6pc}\xymatrixcolsep{0.3pc} \xymatrix{  &1
\ar[dl] \ar[dr] &&5 \ar[dl] \ar[dr] \\ 2 &&3 &&2 } }$
\endignore, belongs to $\pinfdeltamod$, this amounts to a routine check
with the aid of \cite{\Kratwo} for instance.  $\qed$
\endexample

\head 6.  The characteristic words of the simple modules \endhead

A {\it segment\/} of a generalized word $w = (p_{i}^{-1} q_i)_{i
\in \ZZ}$ will be any subword of $w$, where we consider the syllables
$p_i^{-1}$ and $q_i$ as indivisible building blocks of
$w$.  So, alternately expressed, the segments of $w$ are just the connected
components of syllables of $w$.  The segments $q_0 p_1^{-1} q_1 \dots$ and
$\dots p_{-1}^{-1} q_{-1} p_0^{-1}$ will be called the {\it principal
segments\/} of $w$; often, we will refer to the former as the {\it principal
right segment\/} of $w$, and to the latter as the {\it principal left
segment\/} for ease of optical reference, even though the distinction of
sides is not substantive (cf\. the remarks following Theorem 5).

Moreover, it will be convenient to say that a generalized word $w$ has finite
projective dimension if the corresponding string module $\St(w)$ has this
property.

Given a simple module $S \in \lamod$, we will construct the (generalized)
word
$w = w(S)$ postulated in Theorem 5 as the limit of a
sequence of successively growing segments of words of finite projective
dimension.  In doing so, we will observe that this process turns periodic
after the construction of (at most) $4n + 1$ pairs of syllables, where $n$
is the
$K$-dimension of
$\la / J$.  Clearly, segments of words of finite projective dimension may
fail to inherit this property. However, we do have

\proclaim{Observations 10}  Let $w$ be a nontrivial finite word.

{\rm (1)}  If $w$ is a segment of a generalized word of finite projective
dimension, then $w$ is a segment of a {\rm finite\/} word of finite
projective dimension.

{\rm (2)}  If $w$ is primitive and the band module $\Bd(w^r, \phi)$ has
finite projective dimension for some $r$ and $\phi$, then, again, $w$ is a
segment of a finite word of finite projective dimension.
\endproclaim

\demo{Proof}  We prove (2) and leave the similar argument backing part (1) to
the reader.  So let $w = p_1^{-1} q_1 \dots p_m^{-1} q_m$ be primitive with
nontrivial flanking syllables $p_1$ and $q_m$.  As in Proposition 1, we let
$\fp_m$ (resp.
$\fq_1$) be the longest paths in
$\Gamma$ such that $\fp_m p_m$ is a path in $K\Gamma \setminus I$ (resp.,
such that $\fq_1 q_1$ is a path in $K \Gamma \setminus I$).  Invoking
Proposition 1, we conclude that the word with graph

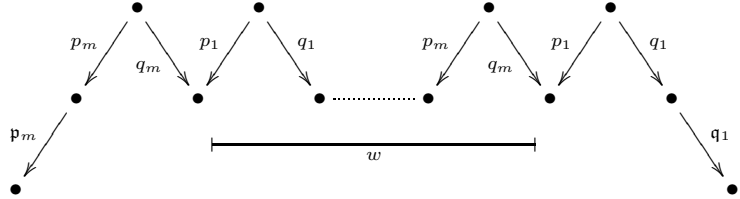
\ignore
$$\xymatrixrowsep{2pc}\xymatrixcolsep{1pc}
\xymatrix{
 &&\bullet \ar[dl]_{p_m} \ar[dr]_(0.55){q_m} &&\bullet \ar[dl]_{p_1}
\ar[dr]^{q_1} &&&&\bullet \ar[dl]_{p_m} \ar[dr]_(0.55){q_m} &&\bullet
\ar[dl]_{p_1}
\ar[dr]^{q_1} \\
 &\bullet \ar[dl]_{\fp_m} &&\bullet \ar@{|-|}[rrrrrr]<-4ex>_{w} &&\bullet
\ar@{.}[rr] &&\bullet &&\bullet &&\bullet
\ar[dr]^{\fq_1}\\
\bullet &&&&& && &&&&&\bullet }$$
\endignore

\noindent has finite projective dimension. \qed
\enddemo

In view of Observations 10, we will  not run any risk of ambiguity if we
henceforth simply refer to `segments of words of finite projective
dimension'.

In order to recognize the algorithmic nature of our construction, the reader
should be familiar with the computation of  projective dimensions of modules
with tree graphs over monomial relation algebras (see \cite{\pre}).  We
precede the description of the procedure with an easy auxiliary statement
spelling out the mechanism of the individual steps.  In view of Proposition
1, the proof is straight-forward and will be omitted.

The notationally somewhat involved second parts of statements (A) and (B)
below, under the heading `more detail', are only relevant for algorithmic
purposes and do not impinge on the further development of the theory; the
reader only interested in the latter is advised to skip
them.  Recall that, for any path $p$ of
positive length in $K\Gamma \setminus I$, we denote by $\fp$ the unique
longest path with the property that the concatenation $\fp p$ is still a
path in
$K\Gamma \setminus I$.

\proclaim{Observations 11}

{\rm (A)}  Suppose that $u = q_{-t} p_{-(t-1)}^{-1} \dots q_{t-1} p_t^{-1}$
with $t \ge 1$ is a segment of a centered word of finite projective
dimension; so in particular, $u$ is nontrivial if and only if the segment
$p_0^{-1} q_0$ is nontrivial. 

 Then there exist unique
\underbar{shortest} paths
$q_t$ and
$p_{-t}$ such that $\, (p_{-t})^{-1}uq_t\, $ is again a segment of a word of
finite projective dimension.  In fact, the path $q_t$ depends only on the
first arrow of $p_t$ in case the latter path has positive length, and is
trivial otherwise; symmetrically, $p_{-t}$ depends only on the first arrow
of $q_{-t}$ in case that path has positive length, and is trivial otherwise.

(Note that, in general, the choice $q_t = e(t)$ will be ruled out, since all
syllables to the right of a trivial syllable $q_i$ with
$i \ge 0$ are required to be trivial by the definition of a centered
word.  Analogous considerations apply to the left-hand side.)
\smallskip

More detail: If $p_t = p_t e(t)$ has positive length and $\la e(t)$ has graph

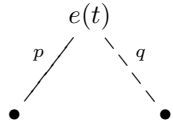
\ignore
$$\xymatrixrowsep{2pc}\xymatrixcolsep{1pc}
\xymatrix{
 &e(t) \edge[dl]_{p} \ar@{--}[dr]^{q} \\
\bullet &&\bullet }$$
\endignore

\noindent where $p$ contains $p_t$ as a right subpath and $q$ is a path of
length $\ge 0$, then $q_t$ can be described as follows: It is the shortest
right subpath of $q$ such that, in writing $q =\fq_t q_t$, we either have

$\bullet$  $\St(p^{-1} q_t) \in \pinflamod$, or else

$\bullet$  $\length(q_t) \ge 1$, and there exists a path $r$ of positive
length with the property that
$q_t r^{-1}$ is a word and
$\St(\frak{r}^{-1}\frak{q}_t)$ has finite projective dimension; here
$\frak r$ and $\frak q$ relate to $r$ and $q$, respectively, as indicated
ahead of the lemma.

The path $p_{-t}$ has an analogous description.
\medskip

{\rm (B)}  Now suppose that $u =  p_{-(t-1)}^{-1} q_{-(t-1)}
\dots p_{t-1}^{-1} q_{t-1}$ with $t \ge 1$ is a segment of a centered word of
finite projective dimension.  

Then there exist unique
\underbar{longest} paths
$p_t$ and $q_{-t}$ such that
$\, q_{-t} u (p_t)^{-1}\,$ is a segment of a word of finite projective
dimension.  In fact, the path $p_t$ depends only on the last arrow of
$q_{t-1}$ in case the latter path has positive length, and is trivial
otherwise; symmetrically, $q_{-t}$ depends only on the last arrow of
$p_{-(t-1)}$ in case the latter path has positive length, and is trivial
otherwise.
\smallskip

More detail: If $q_{t-1} = \tilde{e}(t) q_{t-1}$ has positive length and the
injective hull of $\la\tilde{e}(t) / J\tilde{e}(t)$ has graph

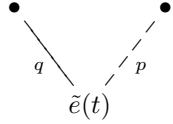
\ignore
$$\xymatrixrowsep{2pc}\xymatrixcolsep{1pc}
\xymatrix{
\bullet \edge[dr]_q &&\bullet \ar@{--}[dl]^p\\
 &\etil(t) }$$
\endignore

\noindent where $q$ contains $q_{t-1}$ as a left subpath, then
$p_t$ can be described as follows:  it is the unique  longest nontrivial
subpath $r$ of $p$ such that

$\bullet$  $\St(\frak{r}^{-1} \frak{q}_{t-1}) \in \pinflamod$, if such a
path $r$ exists,

\noindent and trivial otherwise.

Again, the path $p_{-t}$ can be described analogously.
 \qed
\endproclaim

In contrast to our usual practice, we will consistently record trivial
syllables in the following construction.  

\example{12. Construction of the ``characteristic word'' of a simple module}

Start with a simple left module $S = \la e/ Je$ and construct a generalized
word $w = w(S)$ centered at $e$ as follows:

\smallskip {\bf Step 0:} Choose paths $p_0$ and $q_0$ starting in $e$ and
having minimal length with the property that $w_0 = p_0^{-1} q_0$ is a
segment of a centered word of finite projective dimension. Note that, up to
a swap of roles, $p_0$ and $q_0$ are uniquely determined by this requirement.
\smallskip

{\bf Step t, t $\ge$ 1:}  Suppose that the centered word $w_{t-1} =
p_{-(t-1)}^{-1} q_{-(t-1)} \dots p_{t-1}^{-1} q_{t-1}$ has already been
constructed.   According to part B of Observations 11, we first find the
unique {\it longest} paths $q_{-t}$ and
$p_t$ with the property that $q_{-t} w_{t-1} p_t^{-1}$ is a segment of a
word of finite projective dimension and, according to part A of Observations
11, we then choose $p_{-t}$ and $q_t$ as the unique {\it shortest} paths
such that $w_t = p_{-t}^{-1} q_{-t} w_{t-1} p_t^{-1} q_t$ is a segment of a
word of finite projective dimension.

\smallskip

In Step  $2n$, at the latest, we hit periodicity on both sides. We explain
this for the principal right segment, the left-hand side behaving
symmetrically.  If $q_{2n}$ has length zero, our claim is trivial, since in
that case all further syllables on the right are paths of length zero.  So
suppose that $\length(q_{2n}) \ge 1$. Then the principal right segment of
$w_{2n}$, that is, the word
$q_0 p_1^{-1} q_1 \dots p_{2n}^{-1} q_{2n}$ consists of $4n+1$ nontrivial
syllables and thus gives rise to a string module with socle dimension
$2n+1$, meaning that some simple, say $\la
\tilde{e} / J \tilde{e}$, occurs with multiplicity at least $3$ in this
socle.  Hence at least two of the terminal arrows of the corresponding
nontrivial paths $q_i$ ending in $\tilde{e}$ coincide, say $q_k$ and $q_l$
with $k < l \le 2n$.  Thus Observations 11 guarantee that $p_{k+1} =
p_{l+1}$, $q_{k+1} = q_{l+1}$, $p_{k+2} = p_{l+2}$, and so forth.

This completes the description of the construction. \endexample

Even though the left and right periodic centered word $w = w(S)$ produced by
this construction is only determined up to inversion, we will use the
definite article in referring to it.

\definition{Definition 13}  Given a simple left $\la$-module $S = \la e /
Je$, the centered word
$w = w(S)$ of Construction 12 will be called the {\it characteristic word\/}
of $S$, and the module $\St(w)$ the {\it characteristic phantom\/} of $S$.
\enddefinition

This terminology is justified by

\proclaim{Proposition 14} The characteristic words of the simple left
$\la$-modules are (up to inversion) the centered words $w_1,
\dots, w_n$ postulated in Theorem 5.
\endproclaim

Note that proving Proposition 14 will, at the same time, establish our main
result, Theorem 5.  This will be done in the next section.  As is clear from
their construction, the characteristic words of the simple modules are strongly
interconnected.  We record this fact in

\proclaim{Remark 15}  Let $w = (p_i^{-1} q_i)_{i \in \ZZ}$ be the
characteristic word of a simple left $\la$-module $S$ and $e(i)$ the coinciding
starting point of the paths $p_i$ and $q_i$.  If $i \ge 1$ and $p_{i-1}$ is
nontrivial, then the segment $\, q_i p_{i+1}^{-1} \dots\,$ of $w$ is, up to
re-indexing, one of the principal segments of the characteristic word of $\la
e(i) / J e(i)$.  Analogously, if $i
\le -1$ and $q_{i-1}$ is nontrivial, then $\, \dots q_{i-2}^{-1} p_{i-1}
q_{i-1}^{-1} p_i \, $ is, up to re-indexing, a principal segment of the
characteristic word of $\la e(i) / J e(i)$. \qed
\endproclaim

The paths arising as syllables of the characteristic words play a pivotal
role in the structural makeup of arbitrary modules of finite projective
dimension.  This is reflected by the following facts.

\proclaim{Proposition 16} Let $w = (p_i^{-1} q_i)_{i
\in \ZZ}$ be the characteristic word of $S = \la e / Je$, and again denote
the joint starting point of $p_i$ and $q_i$ by $e(i)$.  Moreover, let $M$ be
any object of
$\pinflamod$.

{\rm (A)}  Suppose that $x$ is a top element of $M$.  If $x$ is of type
$e=e(0)$, then
$p_0 x \ne 0$ and $q_0 x \ne 0$.  If $x$ is of type $e(i)$ for some $i \ge
1$, then $q_ix \ne 0$, and if $x$ is of type $e(i)$ for some $i \le -1$, then
$p_i x \ne 0$.

{\rm (B)}  Now suppose that the path $p_i$ is nontrivial and $x = e(i)x$ for
some $i \ge 1$ is an element of $M$  with the property that
$p_i x$ is a nonzero element of $\soc M \cap q_{i-1} M$.  Then $x$ is a top
element of $M$.  The same conclusion holds if $i \le -1$, the path
$q_i$ is nontrivial, and $q_i x$ is a nonzero element of $\soc M \cap p_{i+1}
M$.
\endproclaim

\demo{Proof}  We may clearly restrict our attention to the case where $M$ is
a string or a band module.  In the former case, our assertions are immediate
consequences of our choices of the $p_i$ and $q_i$.  So suppose that $M$ is
a band module based on a primitive word $v$.  By Observation 10(2), finite
projective dimension of $M$ forces $v$ to be a segment of a word of finite
projective dimension, whence our claims again follow from Observations 11
and Construction 12. \qed \enddemo

We conclude this section with examples demonstrating that all theoretically
possible scenarios actually occur: left and right termination of a
characteristic word $w$, meaning that, for $i \gg 0$, the syllables $p_i$,
$q_i$ and $p_{-i}$, $q_{-i}$ are primitive idempotents (recall that, in view
of Theorem 5 and Proposition 14, this occurs precisely when the corresponding
simple module $S = \la e / Je$ has a $\pinflamod$-approximation); onesided
termination of $w$; non-trivial left and right periodicity,
$$w = \dots uuu  *** \dots ***  vvv
\dots,$$ 
with primitive words $u$ and $v$ which are devoid of common
syllables;
 and periodicity in the strongest possible sense, i\.e\., $w = \dots vvv
\dots$, where $v$ is a primitive word.

\definition{Example 17}  Let $\Gamma$ be the quiver

$$
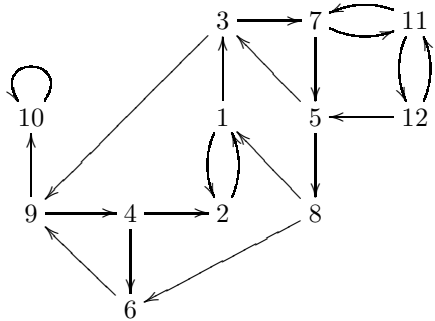
\xy \xymatrix {&&  3\ar[r]\ar[ddll] &7\ar[d]\ar@/_/[r]
&11\ar@/_/[l]\ar@/_/[d]\\ 10
\ar@(ur,ul)&&1\ar[u]\ar@/_/[d]&5\ar[d]\ar[ul]&12\ar[l]\ar@/_/[u]\\
9\ar[r]\ar[u]&4\ar[r]\ar[d]&2\ar@/_/[u]&8\ar[dll]\ar[ul]\\ &6 \ar[ul]   
&&&\\}
\endxy
$$


\noindent and $\la = K\Gamma/I$, where the ideal $I \subseteq K\Gamma$ is
chosen so that the indecomposable projective left $\la$-modules have graphs

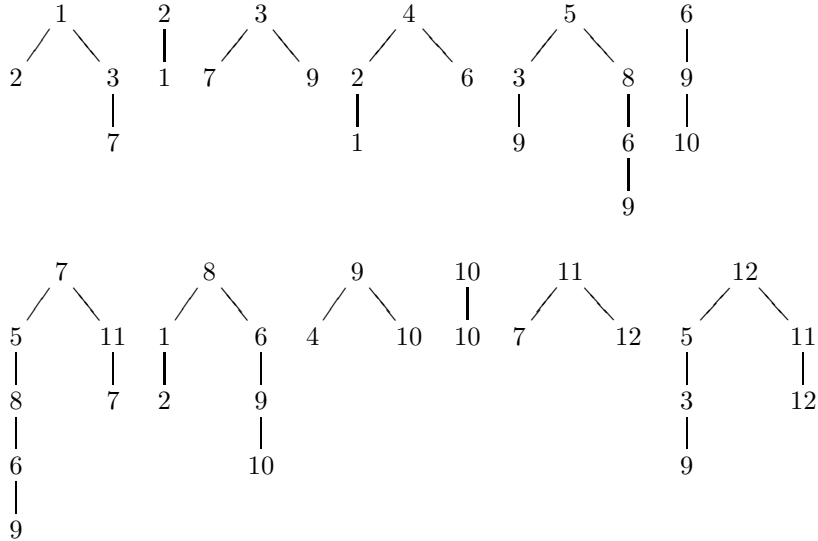
\ignore
$$\xymatrixrowsep{1pc}\xymatrixcolsep{0.5pc}
\xymatrix{
 &1 \edge[dl] \edge[dr] &&2 \edge[d] &&3 \edge[dl] \edge[dr] &&&4
\edge[dl] \edge[dr] &&&5 \edge[dl] \edge[dr] &&6 \edge[d]\\ 2 &&3 \edge[d]
&1 &7 &&9 &2 \edge[d] &&6 &3 \edge[d] &&8 \edge[d] &9
\edge[d]\\
 &&7 &&&&&1 &&&9 &&6 \edge[d] &10\\
 &&&& &&&& &&&&9\\
 &7 \edge[dl] \edge[dr] &&&8 \edge[dl] \edge[dr] &&&9 \edge[dl]
\edge[dr] &&10 \edge[d] &&11 \edge[dl] \edge[dr] &&&12 \edge[dl]
\edge[dr]\\ 5 \edge[d] &&11 \edge[d] &1 \edge[d] &&6 \edge[d] &4 &&10 &10 &7
&&12 &5
\edge[d] &&11 \edge[d]\\ 8 \edge[d] &&7 &2 &&9 \edge[d] &&&& &&&&3 \edge[d]
&&12\\ 6 \edge[d] &&&&&10 &&&& &&&&9\\ 9 }$$
\endignore

\noindent Then $\la$ is a string algebra with simple left modules $S_i =
\la e_i/ Je_i$, $1 \le i \le 12$.  One readily finds that the characteristic
phantoms (in the sense of the above definition  --  see also Theorem 5) of
$S_1$, $S_7$, and $S_8$ are as follows:

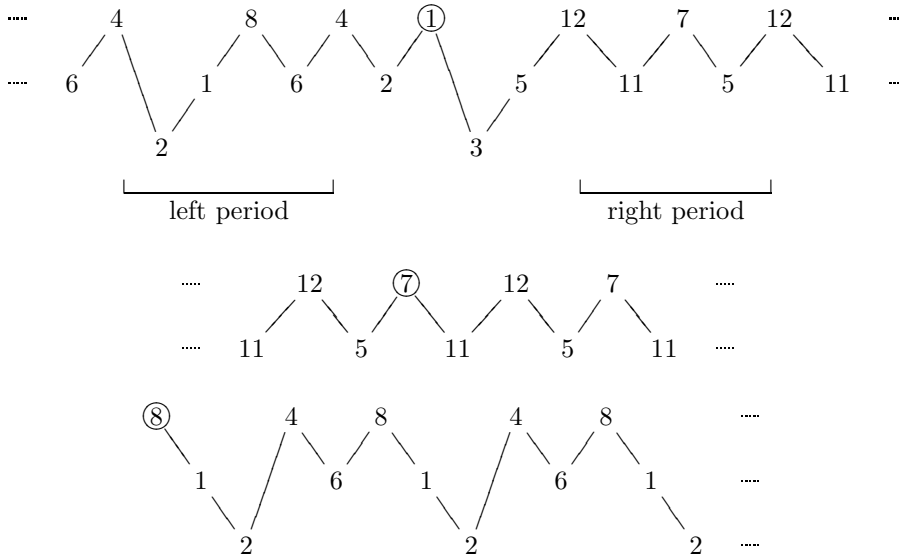
\ignore
$$\xymatrixrowsep{1pc}\xymatrixcolsep{0.5pc}
\xymatrix{
\ar@{.}[r] &&&4 \edge[dl] \edge[ddr] &&&8 \edge[dl] \edge[dr] &&4
\edge[dl] \edge[dr] &&\encirc{1} \edge[dl] \edge[ddr] &&&12 \edge[dl]
\edge[dr] &&7 \edge[dl] \edge[dr] &&12 \edge[dl] \edge[dr] && \ar@{.}[r] &\\
\ar@{.}[r] &&6 &&&1 \edge[dl] &&6 &&2 &&&5 \edge[dl] &&11 &&5 &&11 &
\ar@{.}[r] &\\
 &&& \ar@^{|-|}[rrrrr]<-4ex>_{\txt{left period}}  &2 &&&&&&&3 &&
\ar@^{|-|}[rrrr]<-4ex>_{\txt{right period}} &&&& }$$
\endignore

\ignore
$$\xymatrixrowsep{1pc}\xymatrixcolsep{0.5pc}
\xymatrix{
\ar@{.}[r] &&&12 \edge[dl] \edge[dr] &&\encirc{7} \edge[dl] \edge[dr] &&12
\edge[dl] \edge[dr] &&7 \edge[dl] \edge[dr] &&\ar@{.}[r] &\\
\ar@{.}[r] &&11 &&5 &&11 &&5 &&11 & \ar@{.}[r] & }$$
\endignore

\ignore
$$\xymatrixrowsep{1pc}\xymatrixcolsep{0.5pc}
\xymatrix{ \encirc{8} \edge[dr] &&&4 \edge[ddl] \edge[dr] &&8 \edge[dl]
\edge[dr] &&&4 \edge[ddl] \edge[dr] &&8 \edge[dl] \edge[dr] &&& \ar@{.}[r]
&\\
 &1 \edge[dr] &&&6 &&1 \edge[dr] &&&6 &&1 \edge[dr] && \ar@{.}[r] &\\
 &&2 &&&&&2 &&&&&2 & \ar@{.}[r] & }$$
\endignore

\noindent In each case, the center is circled.

In particular, we see that the characteristic word $w_1$ of $S_1$ is
two-sided infinite with distinct left and right periods, the characteristic
word $w_7$ of
$S_7$ is of the form $w_7 = \dots vvv \dots$, where $v$ is a primitive word,
while the characteristic word
$w_8$ of $S_8$ is infinite only on one side.  Note moreover that the
characteristic word of $S_7$ coincides with that of $S_{12}$, while that of
$S_4$ results from that of $S_8$ through deletion of the first two syllables.

In view of Theorem 5, $\pinflamod$ fails to be contravariantly finite in
$\lamod$.  The theorem, in fact, supplies more precise information:
  Namely, the simple modules $S_2$, $S_3$, $S_6$, $S_9$, $S_{10}$, and
$S_{11}$ are precisely those having $\pinflamod$-approximations and, for
each of the listed indices $i$, the minimal
$\pinflamod$-approximation of $S_i$ coincides with $\la e_i$. \qed
\enddefinition

\definition{Example 18}  Let $\Gamma'$ be the quiver with $12$ vertices
resulting from that of Example 17 through deletion of two arrows, those from
$8$ to
$1$ and from $12$ to $11$.  We define the string algebra $\la' = K\Gamma' /
I'$, where
$I'
\subseteq K \Gamma'$ is chosen in such a way that the graphs of the
indecomposable projective left $\la'$-modules $\la' e_i$, for $i \in \{1,
\dots, 12\}
\setminus\{8,12\}$ coincide with those of the corresponding $\la e_i$-modules
of Example 17, while $\la' e_8$ and $\la' e_{12}$ have graphs

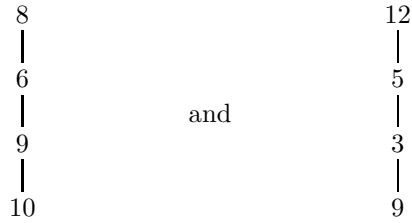
\ignore
$$\xymatrixrowsep{1pc}\xymatrixcolsep{5pc}
\xymatrix{ 8 \edge[d] &&12 \edge[d]\\ 6 \edge[d] &\dropdown{3}{\txt{and}} &5
\edge[d]\\ 9 \edge[d] &&3 \edge[d]\\ 10 &&9 }$$
\endignore

\noindent respectively.  Then the characteristic phantoms of the simple left
$\la'$-modules are as follows, the centers being again highlighted:

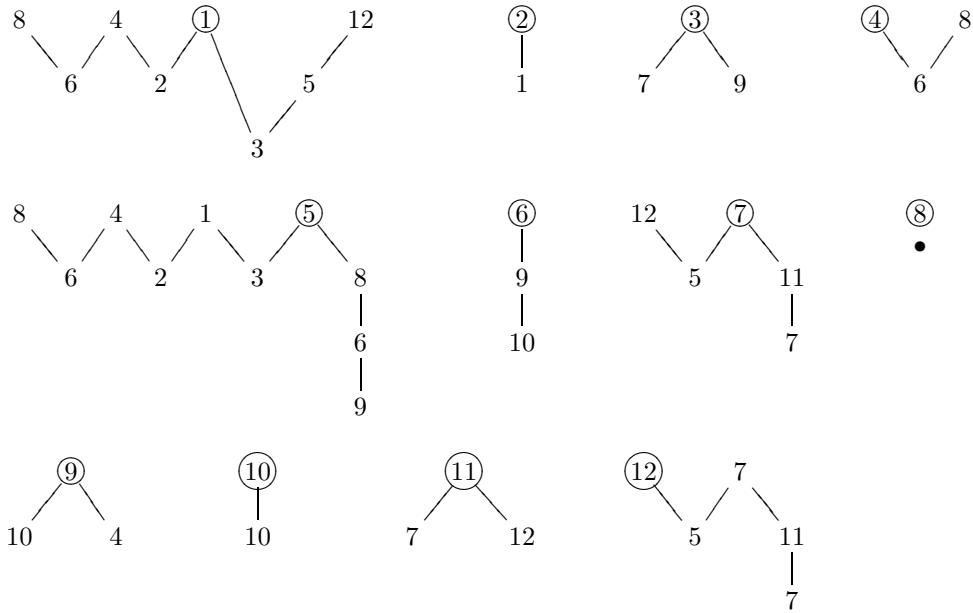
\ignore
$$\xymatrixrowsep{1pc}\xymatrixcolsep{0.5pc}
\xymatrix{ 8 \edge[dr] &&4 \edge[dl] \edge[dr] &&\encirc{1} \edge[dl]
\edge[ddr] &&&12 \edge[dl] &&&\encirc{2} \edge[d] &&&&\encirc{3} \edge[dl]
\edge[dr] &&&&\encirc{4} \edge[dr] &&8 \edge[dl]\\
 &6 &&2 &&&5 \edge[dl] &&&&1 &&&7 &&9 &&&&6\\
 &&&&&3\\ 8 \edge[dr] &&4 \edge[dl] \edge[dr] &&1 \edge[dl] \edge[dr]
&&\encirc{5}
\edge[dl] \edge[dr] &&&&\encirc{6} \edge[d] &&&12 \edge[dr] &&\encirc{7}
\edge[dl] \edge[dr] &&&&\encirc{8} \dropdown{3}{\bullet}\\
 &6 &&2 &&3 &&8 \edge[d] &&&9 \edge[d] &&&&5 &&11 \edge[d]\\
 &&&&&&&6 \edge[d] &&&10 &&&&&&7\\
 &&&&&&&9\\
 &\encirc{9} \edge[dl] \edge[dr] &&&&\encirc{{10}} \edge[d] &&&&\encirc{{11}}
\edge[dl] \edge[dr] &&&&\encirc{{12}} \edge[dr] &&7 \edge[dl]
\edge[dr] \\  10 &&4 &&&10 &&&7 &&12 &&&&5 &&11 \edge[d]\\
 &&&&& &&&&& &&&&& &7 }$$
\endignore

\noindent Thus Theorem 5 tells us that $\pinflamod$ is contravariantly
finite in $\lamod$, and that the displayed characteristic phantoms are the
minimal
$\pinflamod$-approximations of the simple modules corresponding to the
centers. \qed
\enddefinition

\head{7. Proof of the main result} \endhead

Our plan is to establish Theorem 5 by way of proving Proposition 14.

Throughout this section, we fix a simple left module $S = \la e/ Je$ with $e
\in \{e_1, \dots, e_n\}$ and let $w = w(S) = (p_i^{-1}q_i)_{i \in \ZZ}$ be
the characteristic word of $S$ as described in Construction 12 of Section
6.  To cope with the ambiguity arising from the fact that characteristic
words are only unique up to inversion, we will, in the sequel, refer to the
orientation of this fixed choice of $w$ as {\it normalized\/}.
 Since all parts of the theorem involving the word
$w$ are true if the latter is trivial (recall that this occurs precisely when
$S$ has finite projective dimension), we will henceforth assume that $w$ is
nontrivial, i\.e., that at least one of the paths $p_0$, $q_0$ is
nontrivial.  Relative to the standardized sequence of top elements $(x_i)_{i
\in \supp(w)}$ of $\St(w)$, as introduced in Section 2, the string module
$\St(w)$ has graph

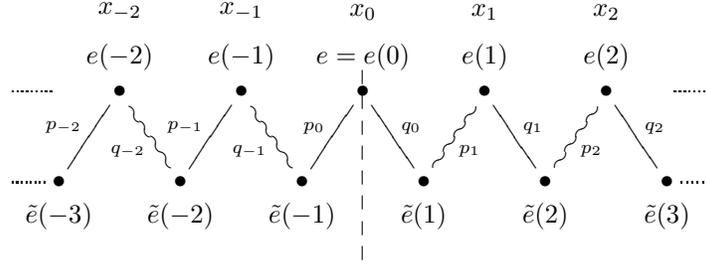
\ignore
$$\xymatrixrowsep{2pc}\xymatrixcolsep{1pc}
\xymatrix{
\ar@{.}[r] &&\bullet \dropup{3}{e(-2)} \dropup{7}{x_{-2}} \edge[dl]_{p_{-2}}
\wiggledge[dr]_(0.55){q_{-2}} &&\bullet
\dropup{3}{e(-1)} \dropup{7}{x_{-1}}
\edge[dl]_{p_{-1}}
\wiggledge[dr]_(0.55){q_{-1}} &&\bullet \save+<0ex,3ex> \drop{e=e(0)}
\ar@{--}[dd] \restore
\dropup{7}{x_{0}} \edge[dl]_{p_0}
\edge[dr]^{q_0} &&\bullet  \dropup{3}{e(1)} \dropup{7}{x_1}
\wiggledge[dl]^(0.55){p_1} \edge[dr]^{q_1} &&\bullet \dropup{3}{e(2)}
\dropup{7}{x_2}
\wiggledge[dl]^(0.55){p_2}
\edge[dr]^{q_2} &\ar@{.}[r] &\\
\ar@{.}[r] &\bullet \dropdown{3}{\etil(-3)} &&\bullet
\dropdown{3}{\etil(-2)} &&\bullet
\dropdown{3}{\etil(-1)} &&\bullet \dropdown{3}{\etil(1)} &&\bullet
\dropdown{3}{\etil(2)} &&\bullet \dropdown{3}{\etil(3)}
\ar@{.}[r] &\\
 &&&&&&
 }$$
\endignore

\noindent Here we denote by $e(i)$ the starting point of the path
$q_i$ for $i \ge 0$, and by $\tilde{e}(i+1)$ its end point;  then $e(i)$ and
$\tilde{e}(i)$ are also the starting and end points of the
$p_i$ for $i \ge 1$, respectively; moreover, $e(0) = e$.  The principal left
segment of $w$ is labeled similarly, as shown in the graph of
$\St(w)$.  Contrasting our usual convention, we have marked some of the
edges by wiggly, as opposed to straight, lines to emphasize the following
difference in roles: the straight edges indicate paths chosen as short as
possible without forfeiting finite projective dimension of
$\St(w)$, whereas the paths represented by wiggly edges are chosen as long
as possible under this restriction.  As in the previous sections,
$\varphi: \St(w) \rightarrow S$ denotes the canonical map which sends
$x_0$ to $e + Je$, while sending the other $x_i$ to zero.

We smooth the road towards a proof of Theorem 5 with a final definition and
a few auxiliary facts.

\definition{Definition 19}  Two nontrivial centered words
$\what = (\phat_i^{-1}\qhat_i)_{i \in \ZZ}$ and $\wtil =
(\ptil_i^{-1}\qtil_i)_{i \in \ZZ}$ are said to {\it have the same
orientation} if they are centered in the same primitive idempotent and
either $\phat_0$ and $\ptil_0$ have the same first arrow, or else $\qhat_0$
and $\qtil_0$ have the same first arrow. (Note that, if all of the paths
$\phat_0$, $\ptil_0$, $\qhat_0$,
$\qtil_0$ are nontrivial, the condition that $\what$ and $\wtil$ have the
same orientation is equivalent to the requirement that
$\phat_0$ and $\ptil_0$, as well as $\qhat_0$ and $\qtil_0$, share first
arrows.)

Let $S$ and $w = (p_i^{-1}q_i)$ be as fixed at the beginning of this
section. If $\what = (\phat_i^{-1}\qhat_i)$ has the same orientation as $w$,
we call $\hat{q}_i$ (resp. $\hat{p}_i$) a {\it right discontinuity} of
$\hat{w}$ relative to $w$ in case $i \ge 0$ and $\hat{q}_i \ne q_i$ (resp.,
$i \ge 1$ and $\hat{p}_i \ne p_i$); {\it left discontinuities} are defined
symmetrically. Both types are also briefly called discontinuities of $\what$.
\enddefinition

In the sequel, all discontinuities will be relative to the fixed
characteristic word $w$ of $S = \la e / Je$.

\proclaim{Lemma 20}  Let $\hat{v} =\hat{p}_0^{-1}\hat{q}_0 \dots
\hat{p}_t^{-1}\hat{q}_t$ be a primitive word with $t \ge 0$ and all of the
listed syllables nontrivial; moreover, suppose that the joint starting point
of $\phat$ and $\qhat$ equals $e$.  Expand
$\hat{v}$ to a twosided infinite word $\hat{w} = \dots
\hat{v}\hat{v}\hat{v} \dots$ centered in $e$ as illustrated by the following
graph

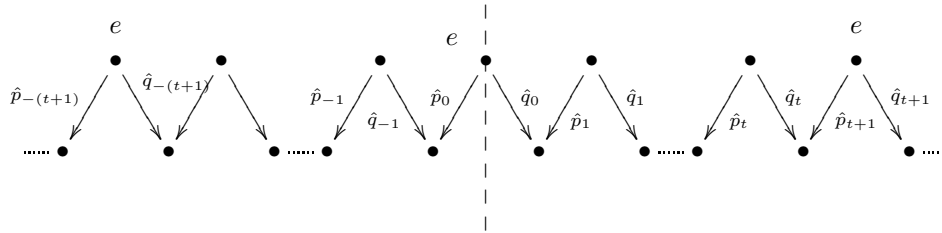
\ignore
$$\xymatrixrowsep{2pc}\xymatrixcolsep{0.75pc}
\xymatrix{
 &&\bullet \dropup{3}{e} \ar[dl]_{\hatp_{-(t+1)}}
\ar[dr]^(0.35){\hatq_{-(t+1)}} &&\bullet
\ar[dl] \ar[dr] &&&\bullet \ar[dl]_{\hatp_{-1}} \ar[dr]_(0.55){\hatq_{-1}}
&&\bullet
\save+<-3ex,2ex> \drop{e} \restore \save+<0ex,5ex> \ar@{--}[dd] \restore
\ar[dl]_{\hatp_0} \ar[dr]^{\hatq_0} &&\bullet \ar[dl]^(0.55){\hatp_1}
\ar[dr]^{\hatq_1} &&&\bullet \ar[dl]^(0.55){\hatp_t} \ar[dr]^{\hatq_t}
&&\bullet \dropup{3}{e}
\ar[dl]^(0.55){\hatp_{t+1}} \ar[dr]^{\hatq_{t+1}}\\
\ar@{.}[r] &\bullet &&\bullet &&\bullet \ar@{.}[r] &\bullet &&\bullet
&&\bullet &&\bullet
\ar@{.}[r] &\bullet &&\bullet &&\bullet \ar@{.}[r] &\\
 &&&&&&&&& }$$
\endignore

\noindent where $\phat_i = \phat_j$ and $\qhat_i = \qhat_j$ whenever $i$ is
congruent to $j$ modulo $t+1$.  Finally, suppose that $\what$ has the same
orientation as $w$.

If $\what$ has a right discontinuity, then the first right discontinuity of
$\what$ is among the paths $\qhat_0, \dots,
\qhat_t$, $\phat_1, \dots, \phat_{t+1}$.  Similarly, if $\what$ has a left
discontinuity, then the first such is among $\phat_0,
\dots, \phat_{-t}$, $\qhat_{-1}, \dots, \qhat_{-(t+1)}$.

Now suppose that $\what$ is a word of finite projective dimension.  If $\what$ has
a right discontinuity, and the first right discontinuity is $\qhat_i$ for some
$i$, then the path $\qhat_i$ contains the path $q_i$  as a proper right
subpath; if, on the other hand, the first right discontinuity of $\what$ is
$\phat_i$, the latter path is contained in the path $p_i$ as a proper left
subpath.  In case of existence, the first left discontinuity of $what$ is subject
to mirror-symmetric conditions.
\endproclaim

\demo{Proof}  To verify the first assertion, suppose that
$\qhat_i = q_i$ for $0 \le i \le t$ and $\phat_i = p_i$ for $1 \le i \le
t+1$.  Then the paths $q_0$ and $p_{t+1}$ starting in $e$ are both
nontrivial, and the first arrow of $p_{t+1}$ differs from the first arrow of
$q_0 = \qhat_0 = \qhat_{t+1}$.  By Construction 12, this implies that
$q_{t+1} = q_0$, i\.e., $\qhat_{t+1} = q_{t+1}$, and further that $p_{t+2} =
p_1 = \phat_1 = \phat_{t+2}$,
$q_{t+2} = q_1 = \qhat_1 = \qhat_{t+2}$, and so forth, meaning that the
principal right segments of $w$ and $\what$ coincide.  The counterpart
dealing with the principal left segment of $\what$ is symmetric.

The final assertions, concerning first right and left discontinuities in
case $\what$ is a word of finite projective dimension, are immediate
consequence of the construction of
$w$ (Section 6, Construction 12).
\qed
\enddemo

In the proof of the next lemma, it will turn out handy that every
homomorphism from a pseudo-band module $\Bd(v^r, \phi)$ to $S$ factors
through an `expanded' pseudo-band module $\Bd(v^{2r},
\psi)$.

\proclaim{Remark 21}  Since $K$ is an infinite field, any pseudo-band module
$\Bd(v^r, \phi)$ is contained as a direct summand in a pseudo-band module
$\Bd(v^{s}, \psi)$, for any integer $s \ge r$.
\endproclaim

\proclaim{Lemma 22}  Let $\vhat$ and $\what$ be as in the blanket hypothesis
of Lemma 20 and retain all of the notation introduced there.  Moreover,
suppose that $\what$ is a word of finite projective dimension having both
right and left discontinuities.

Then, given any pseudo-band module $B =
\Bd(\vhat^r, \phi)$, say with graph

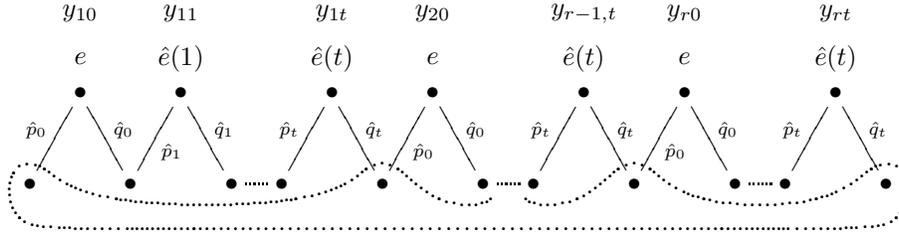
\ignore
$$\xymatrixrowsep{2pc}\xymatrixcolsep{0.67pc}
\xymatrix{
 &\bullet \dropup{3}{e} \dropup{7}{y_{10}} \edge[dl]_{\hatp_0}
\edge[dr]^{\hatq_0} &&\bullet
\dropup{3}{\ehat(1)} \dropup{7}{y_{11}} \edge[dl]^(0.55){\hatp_1}
\edge[dr]^{\hatq_1} &&&\bullet \dropup{3}{\ehat(t)} \dropup{7}{y_{1t}}
\edge[dl]_{\hatp_t} \edge[dr]^{\hatq_t} &&\bullet \dropup{3}{e}
\dropup{7}{y_{20}} \edge[dl]^(0.55){\hatp_0}
\edge[dr]^{\hatq_0} &&&\bullet \dropup{3}{\ehat(t)} \dropup{7}{y_{r-1,t}}
\edge[dl]_{\hatp_t} \edge[dr]^{\hatq_t} &&\bullet \dropup{3}{e}
\dropup{7}{y_{r0}}
\edge[dl]^(0.55){\hatp_0}
\edge[dr]^{\hatq_0} &&&\bullet  \dropup{3}{\ehat(t)} \dropup{7}{y_{rt}}
\edge[dl]_{\hatp_t} \edge[dr]^{\hatq_t}\\
\bullet &&\bullet &&\bullet  \ar@{.}[r] &\bullet &&\bullet &&\bullet
\ar@{.}[r] &\bullet
\save[0,0]+(-1,-2);[0,-1]+(1,-2) **\crv{~*=<2.5pt>{.} [0,0]+(0,-3)
&[0,1]+(0,-3) &[0,2]+(-2,3) &[0,2]+(2,3) &[0,3]+(0,-3) &[0,6]+(0,-3)
&[0,7]+(-2,3) &[0,7]+(3,3) &[0,7]+(3,-6) &[0,6]+(0,-6) &[0,-9]+(0,-6)
&[0,-10]+(-3,-6) &[0,-10]+(-3,3) &[0,-10]+(2,3) &[0,-9]+(0,-3)
&[0,-4]+(0,-3) &[0,-3]+(-2,3) &[0,-3]+(2,3) &[0,-2]+(0,-3) &[0,-1]+(0,-3)
}\restore
 &&\bullet &&\bullet  \ar@{.}[r] &\bullet &&\bullet }$$
\endignore
\vskip0.3truein

\noindent relative to a standardized sequence $y_{10}, \dots, y_{rt}$ of top
elements, the homomorphism $f: B \rightarrow S$ sending $y_{10}$ to $e + Je$
and the other $y_{ij}$ to zero factors through the canonical map $\varphi:
\St(w) \rightarrow S$.
\endproclaim

\demo{Proof}  We start by formalizing the graphical information provided:
$\hat{q}_j y_{ij} =
\hat{p}_{j+1} y_{i,j+1}$ for
$j < t$,
$\qhat_t y_{it} = \phat_{i+1,0} y_{i+1,0}$ when $i < r$, and $\qhat_t y_{rt}
=
\sum_{i = 1}^r c_i \phat_0 y_{i0}$, where
$$\pmatrix 0&\cdots&0&c_1\\ 1&\ddots &&\vdots\\ &\ddots&0&\vdots\\
0&\cdots&1&c_r \endpmatrix$$ is the Frobenius companion matrix of the cyclic
automorphism $\phi$.

To facilitate visualization, we once more give the graph of the word
$\what = \dots \vhat \vhat \vhat \dots$,

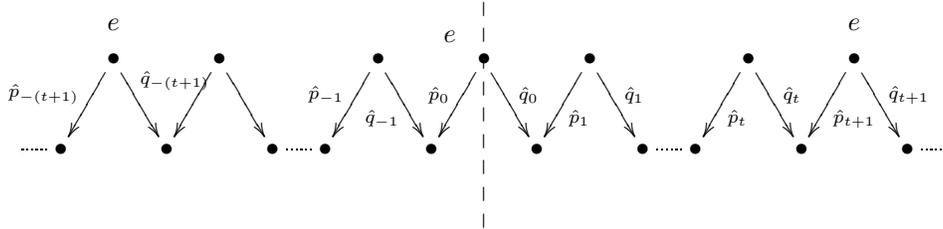
\ignore
$$\xymatrixrowsep{2pc}\xymatrixcolsep{0.75pc}
\xymatrix{
 &&\bullet \dropup{3}{e} \ar[dl]_{\hatp_{-(t+1)}}
\ar[dr]^(0.35){\hatq_{-(t+1)}} &&\bullet
\ar[dl] \ar[dr] &&&\bullet \ar[dl]_{\hatp_{-1}} \ar[dr]_(0.55){\hatq_{-1}}
&&\bullet
\save+<-3ex,2ex> \drop{e} \restore \save+<0ex,5ex> \ar@{--}[dd] \restore
\ar[dl]_{\hatp_0} \ar[dr]^{\hatq_0} &&\bullet \ar[dl]^(0.55){\hatp_1}
\ar[dr]^{\hatq_1} &&&\bullet \ar[dl]^(0.55){\hatp_t} \ar[dr]^{\hatq_t}
&&\bullet \dropup{3}{e}
\ar[dl]^(0.55){\hatp_{t+1}} \ar[dr]^{\hatq_{t+1}}\\
\ar@{.}[r] &\bullet &&\bullet &&\bullet \ar@{.}[r] &\bullet &&\bullet
&&\bullet &&\bullet
\ar@{.}[r] &\bullet &&\bullet &&\bullet \ar@{.}[r] &\\
 &&&&&&&&& }$$
\endignore

\noindent where again $\phat_i = \phat_j$ and $\qhat_i = \qhat_j$ whenever
$i \equiv j \pmod{t+1}$.  The word $\what$ having the same orientation as
$w$, Lemma 20 guarantees that the first right discontinuity of $\what$ is
among the paths $q_i$, $p_j$ with $0
\le i \le t$ and $1 \le j \le t+1$, and the first left discontinuity of
$\what$ is among the paths $p_i$, $q_j$ with $0
\ge i \ge -t$ and $1 \ge j \ge -(t+1)$. In factoring the homomorphism $f$
through $\varphi$, we will separately deal with the cases where neither of
the paths $\phat_0$, $\qhat_0$ is a discontinuity of $\what$, where one of
them is, and where both of them are discontinuities.

The last-mentioned case is immediate:  Namely, we  define $g \in
\Hom_{\la}(B,
\St(w))$ by setting $g(y_{10}) = x_0$ and
$g(y_{ij}) = 0$ for $(i,j) \ne (1,0)$.  This is legitimate, since, by Lemma
20, the paths $\phat_0$ and
$\qhat_0$ contain $p_0$ and $q_0$ as proper right subpaths, respectively;
hence $\phat_0 x_0 = \qhat_0 x_0 = 0$, where $(x_i)_{i \in \ZZ}$ is the
standardized sequence of top elements of $\St(w)$ displayed at the beginning
of the section.
\smallskip

\noindent {\bf Case A.}  One of $\phat_0$, $\qhat_0$ is a discontinuity of
$\what$, but not the other; for symmetry reasons, it is harmless to assume
that $\phat_0$ is a (left) discontinuity.  Lemma 20 tells us that $\phat_0$
contains $p_0$ as a proper right subpath, and again we infer
$\phat_0 x_0 = 0$.

{\it Subcase\/}A.1.  $\, \phat_{t+1}$ is the first right discontinuity of
$\what$.

In this case the segment of the word $w$ relevant to our construction has
the form

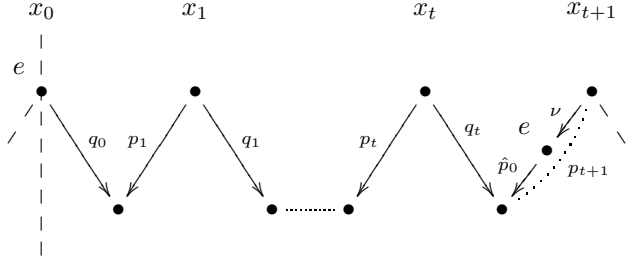
\ignore
$$\xymatrixrowsep{1pc}\xymatrixcolsep{0.5pc}
\xymatrix{
 &\bullet  \save+<-2ex,2ex> \drop{e} \restore \save+<0ex,5ex> \ar@{--}[ddd]
\restore
\dropup{7}{x_0} \ar@{--}[dl]
\ar[ddrr]^{q_0} &&&&\bullet \dropup{7}{x_1} \ar[ddll]_{p_1}
\ar[ddrr]^{q_1} &&&&&&\bullet
\dropup{7}{x_t} \ar[ddll]_{p_t} \ar[ddrr]^(0.4){q_t} &&&&\bullet
\dropup{7} {x_{t+1}} \ar[dl]_{\nu} \ar@{--}[dr]
\ar@{.}@/^/[ddll]^{p_{t+1}} \\
 &&&&&&& &&&&&&&\bullet \save+<-2ex,2ex> \drop{e} \restore
\ar[dl]_(0.45){\hatp_0} && \\
 &&&\bullet &&&&\bullet \ar@{.}[rr] &&\bullet &&&&\bullet \\
 & }$$
\endignore

\noindent where $p_{t+1}  = \phat_{t+1} \nu  = \phat_0 \nu$ for a nontrivial
path
$\nu$, $q_i = \qhat_i$ for $0 \le i \le t$ and $p_i =
\phat_i$ for $1 \le i \le t$.  By Remark 21, we may assume that
$r \ge 2$.  Define
$g \in \Hom_{\la}(B,\St(w))$ as follows, keeping in mind that $c_1 \ne 0$
because
$\phi$ is an automorphism of
$K^r$: Namely, let $g(y_{10}) = x_0 - (c_2/c_1) \nu x_{t+1}$, $g(y_{1j}) =
x_j$ for
$1 \le j
\le t$, $g(y_{20}) = \nu x_{t+1}$, and $g(y_{ij}) = 0$ for $2 \le i \le r$
and
$1 \le j \le t$.  Indeed, in view of the equality $\qhat_0 \nu x_{t+1} = 0$,
it is routine to check that $\qhat_j g(y_{ij}) =
\phat_{j+1} g(y_{i,j+1})$ for $i \ge 1$ and $j \le t-1$, $\qhat_t g(y_{it}) =
\phat_0 g(y_{i+1,0})$ for $i \le r-1$, and $\qhat_t g(y_{rt}) = \sum_{i=1}^r
c_i \phat_0 g(y_{i0})$ under these assignments.  Moreover, we clearly have
$\varphi g = f$.

{\it Subcase\/} A.2.  The first right discontinuity of $\what$ is
$\phat_l$ or $\qhat_l$ with $1 \le l \le t$.

In that case, we have $p_l = \phat_l \nu$, where $\nu$ is a path of length
$\ge 0$, and $\qhat_l = \sigma q_l$ with
$\length(\sigma) \ge 0$, where either $\length(\nu) > 0$ or
$\length(\sigma) > 0$, by Lemma 20.  In either case we obtain
$\qhat_l \nu x_l = 0$, and the following assignments give rise to a
well-defined homomorphism $g \in \Hom_{\la}(B,\St(w))$ with
$\varphi g = f$.  Namely, $g(y_{1j}) = x_j$ for $0 \le j < l$,
$g(y_{1l}) = \nu x_l$, $g(y_{1j}) = 0$ for $l+1 \le j \le t$ (this range
being empty if $l = t$), and $g(y_{ij}) = 0$ for $i \ge 2$ and all $j$.
\smallskip

\noindent {\bf Case B.}  Neither $\phat_0$ nor $\qhat_0$ is a discontinuity
of $\what$.  Let $k \in \{0, \dots t\}$ be such that the first left
discontinuity of $\what$ is either $\phat_{-(k+1)}$ or
$\qhat_{-(k+1)}$, and $l \in \{1, \dots, t+1\}$ such that the first right
discontinuity of $\what$ is either
$\phat_l$ or $\qhat_l$.  Lemma 20 guarantees that either $\phat_{-(k+1)} =
\phat_{t-k}$ contains
$p_{-(k+1)}$ as a proper right subpath, or else $\qhat_{-(k+1)} =
\qhat_{t-k}$ is a proper left subpath of $q_{-(k+1)}$.  In either case we
can write the latter path in the form $q_{-(k+1)} = \qhat_{t-k} \mu$, where
$\mu$ is a path of length $\ge 0$ satisfying the equalities $\phat_{-(k+1)}
\mu x_{-(k+1)} = \phat_{t-k} \mu x_{-(k+1)} = 0$.  Analogously,
$p_l =  \phat_l \nu$, where $\nu$ is a path of length $\ge 0$ such that
$\qhat_l \nu x_l = 0$.  The segment of the characteristic word $w$ which is
decisive for our construction is depicted in the following diagram:

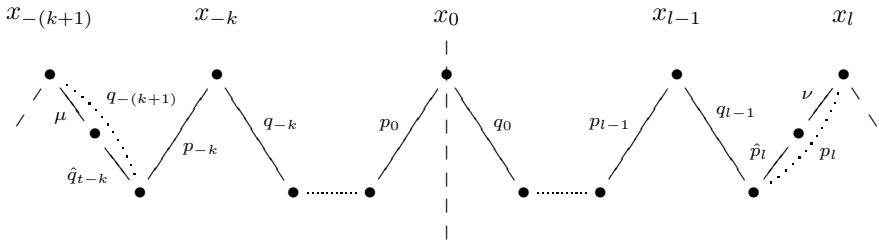
\ignore
$$\xymatrixrowsep{1pc}\xymatrixcolsep{0.5pc}
\xymatrix{
 &\bullet \dropup{5}{x_{-(k+1)}} \ar@{--}[dl] \edge[dr]_(0.55){\mu}
\ar@{.}@/^/[ddrr]^(0.4){q_{-(k+1)}} &&&&\bullet
\dropup{5}{x_{-k}} \edge[ddll]^(0.55){p_{-k}} \edge[ddrr]^{q_{-k}}
&&&&&&\bullet
\dropup{5}{x_0} \save+<0ex,3ex> \ar@{--}[ddd] \restore \edge[ddll]_{p_0}
\edge[ddrr]^{q_0} &&&&&&\bullet \dropup{5}{x_{l-1}} \edge[ddll]_{p_{l-1}}
\edge[ddrr]^(0.4){q_{l-1}} &&&&\bullet \dropup{5}{x_l} \edge[dl]_{\nu}
\ar@{--}[dr] \ar@{.}@/^/[ddll]^{p_l} \\
 &&\bullet \edge[dr]_{\hatq_{t-k}} &&&&& &&&&&& &&&&&&&\bullet
\edge[dl]_{\hatp_l} && \\
 &&&\bullet &&&&\bullet \ar@{.}[rr] &&\bullet &&&&\bullet  \ar@{.}[rr]
&&\bullet  &&&&\bullet
\\
 &&&&&&& &&&& }$$
\endignore

\noindent Here $p_i = \phat_i$ for $-k \le i \le l-1$ and $q_i =
\qhat_i$ for $-k \le i \le l-1$.  Moreover, keep in mind that
$\phat_i = \phat_j$ and $\qhat_i = \qhat_j$ whenever $i \equiv j
\pmod{t+1}$. Remark 21 permits us to assume $r \ge 3$. Referring to the
above graph of $B$ and distinguishing between the cases where $l\le t$ and
$l=t+1$, we can thus define a map $g:
\Bd(v^r,\phi) \rightarrow \St(w)$ as follows:

If $l \le t$, we set
$g(y_{10}) = x_0$, $g(y_{1j}) = x_j$ for $1 \le j \le l-1$, $g(y_{1l}) = \nu
x_l$, $g(y_{1j}) = 0$ for $l+1 \le j \le t$, $g(y_{ij}) = 0$ for $2 \le i \le
r-1$ and all $j$, $g(y_{rj}) = 0$ for $0  \le j \le t-k-1$ (note that this
latter range is empty if $k=t$),
$g(y_{r,t-k}) = c_1
\mu x_{-(k+1)}$, and
$g(y_{rj}) =c_1 x_{-t+j-1}$ for $t-k+1 \le j \le t$ (this range being empty
for $k = 0$).  The verification of the fact that $g$ extends to a
well-defined homomorphism
$B \rightarrow \St(w)$ is a bit tedious  --  various possibilities for $k$
need to be considered separately  --  and we leave it to the reader to fill
in the pertinent computations.

In case $l = t+1$, we can either play the situation back to one of the
previously considered cases by relabeling, or else define $g$ directly as
follows: $g(y_{10}) = x_0 - (c_2/c_1) \nu x_{t+1}$,
$g(y_{1j}) = x_j$ for $1 \le j \le t$, and $g(y_{20}) = \nu x_{t+1}$; for
the pairs $(i,j)$ lexicographically larger than
$(2,0)$, we keep the above specifications of $g(y_{ij})$.  Again one checks
that these definitions give rise to a map $g \in
\Hom_{\la}(B,\St(w))$.

Since clearly $\varphi g = f$ in either case, the proof of the lemma is
complete. $\qed$
\enddemo

\proclaim{Lemma 23} Let $t \ge 0$, and suppose that the syllables
$q_0, \dots, q_t, p_1^{-1}, \dots, p_{t+1}^{-1}$ of the characteristic word
$w = w(S)$ are nontrivial, yielding a primitive word  $v = p_{t+1}^{-1} q_0
p_1^{-1} q_1 \dots p_t^{-1} q_t$ with graph

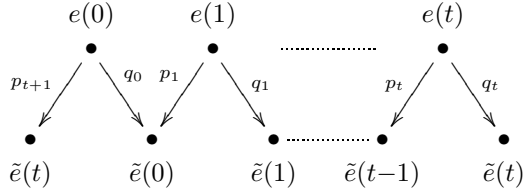
\ignore
$$\xymatrixrowsep{2pc}\xymatrixcolsep{1pc}
\xymatrix{
 &\bullet \dropup{3}{e(0)} \ar[dl]_{p_{t+1}} \ar[dr]^(0.4){q_0} &&\bullet
\dropup{3}{e(1)}
\ar[dl]_(0.4){p_1} \ar[dr]^{q_1} & \ar@{.}[rr] &&&\bullet \dropup{3}{e(t)}
\ar[dl]_{p_t}
\ar[dr]^{q_t} \\
\bullet \dropdown{3}{\etil(t)} &&\bullet \dropdown{3}{\etil(0)} &&\bullet
\dropdown{3}{\etil(1)} \ar@{.}[rr] &&\bullet \dropdown{3}{\etil(t{-}1)}
&&\bullet
\dropdown{3}{\etil(t)}  }$$
\endignore

\noindent  Moreover, suppose that $e(0) = e$ and that the principal right
segment $w_{\text right}$ of $w$ is periodic of the form $w_{\text right} = 
(q_0 p_1^{-1} \dots q_t p_{t+1}^{-1}) (q_0 p_1^{-1} \dots
 q_t p_{t+1}^{-1}) \dots$; in other words,
$$p_{t+1}^{-1} w_{\text right} = vvv \dots.$$

\noindent Then the direct sums
$$\biggl(\bigoplus_{i=0}^t \la e(i)/Je(i)\biggr)^{(\NN)}
\qquad\text{and}\qquad
\biggl( \bigoplus_{i=0}^t \la \etil(i)/J\etil(i) \biggr)^{(\NN)}$$ are
$\pinflamod$-phantoms and $\sinflamod$-phantoms of
$S$.  In particular, $S$ has neither a $\pinflamod$-  nor an
$\sinflamod$-approxi\-ma\-tion.

Moreover, any effective $\pinflamod$- or $\sinflamod$-phantom of $S$ has a
copy of $\biggl(\bigoplus_{i=0}^t \la e(i)/Je(i)\biggr)^{(\NN)}$ in its top
and a copy of $\biggl( \bigoplus_{i=0}^t \la \etil(i)/J\etil(i)
\biggr)^{(\NN)}$ in its socle.
\endproclaim

\demo{Proof}  We will simultaneously verify that the above semisimple
modules are $\pinflamod$- and $\sinflamod$-phantoms of $S$, as all of our
test maps will have sources in $\sinflamod$.

We set
$\tilde{p}_{t+1} =
\frak{p}_{t+1} p_{t+1}$ and
$\tilde{q}_t =
\frak{q}_t q_t$, where $\frak{p}_{t+1}$ and $\frak{q}_t$ are as in
Proposition 1; the path $p_{t+1}$ being nontrivial, this means that
$\frak{p}_{t+1}$ is the longest path such that
$\tilde{p}_{t+1}$ is a path in
$K\Gamma \setminus I$; the path $\tilde{q}_t$ has an analogous description.
Moreover, we consider, for each positive integer $m \ge 2$, the word $u_m =
\tilde{p}_{t+1}^{-1} q_0
\dots p_t^{-1} q_t v^{m-2} p_{t+1}^{-1} q_0 \dots p_t^{-1} \tilde{q}_t$.
Then $u_m$ is a word of finite projective dimension and $\St(u_m)$ has graph

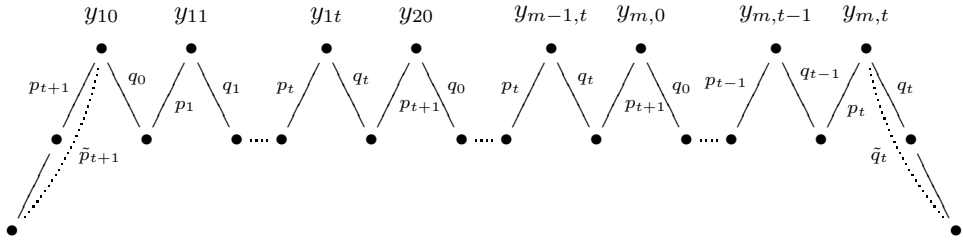
\ignore
$$\xymatrixrowsep{2pc}\xymatrixcolsep{0.5pc}
\xymatrix{
 &&\bullet \dropup{3}{y_{10}} \edge[dl]_{p_{t+1}} \edge[dr]^(0.4){q_0}
\ar@{.}@/^/[ddll]^{\tilde{p}_{t+1}} &&\bullet
\dropup{3}{y_{11}}
\edge[dl]^(0.55){p_1} \edge[dr]^{q_1} &&&\bullet \dropup{3}{y_{1t}}
\edge[dl]_{p_t}
\edge[dr]^(0.4){q_t} &&\bullet
\dropup{3}{y_{20}}
\edge[dl]^(0.55){p_{t+1}} \edge[dr]^{q_0} &&&\bullet  \dropup{3}{y_{m{-}1,t}}
\edge[dl]_{p_t}
\edge[dr]^(0.4){q_t} &&\bullet \dropup{3}{y_{m,0}}
\edge[dl]^(0.55){p_{t+1}} \edge[dr]^{q_0} &&&\bullet \dropup{3}{y_{m,t{-}1}}
\edge[dl]_(0.45){p_{t-1}}
\edge[dr]^(0.35){q_{t-1}} &&\bullet  \dropup{3}{y_{m,t}}
\edge[dl]^(0.6){p_t} \edge[dr]^{q_t} \ar@{.}@/_/[ddrr]_{\tilde{q}_t} \\
 &\bullet \edge[dl] &&\bullet  &&\bullet \ar@{.}[r] &\bullet &&\bullet
&&\bullet \ar@{.}[r] &\bullet &&\bullet  &&\bullet \ar@{.}[r] &\bullet
&&\bullet  &&\bullet \edge[dr] \\
\bullet &&&&& &&&&& &&&&& &&&&& &\bullet }$$
\endignore

\noindent relative to a standardized sequence of top elements $y_{ij}$.  We
equip the set of pairs $(i,j)$ for $1 \le i \le m$ and $0 \le j \le t$ with
the lexicographic order and let
$f: \St(u_m) \rightarrow S$ be the homomorphism with $f(y_{10})) = e + Je$
and
$f(y_{ij}) = 0$ for $(i,j) > (1,0)$.  Next, we will check that any module
$A \in \pinflamod$ with the property that $f$ factors through a map $\rho \in
\Hom_{\la}(A,S)$ contains
$\bigl( \la \etil(j) / J\etil(j) \bigr)^{m-1}$ in its socle and $\bigl( \la
e(j) / Je(j) \bigr)^{m-1}$ in its top, for $j \in \{0, \dots, t\}$.  This
will prove all of our claims.

For that purpose, we let $h \in \Hom_{\la}(\St(u_m), A)$ be such that $\rho
h = f$.  In a first step, we prove, more precisely, that the elements $q_j
h(y_{ij})$, $1 \le i \le m-1$, of $A$ are
$K$-linearly independent for all $j$. Clearly, these elements belong to the
socle of $A$ since the $q_j y_{ij}$'s belong to the socle of $\St(u_m)$ for
$i \le m-1$.  Assume, to the contrary of our claim, that the $q_j
h(y_{ij})$'s, $i \le m-1$, are linearly dependent, and let the pair $(k,l)$
be (lexicographically) minimal with the property that there exists an
equality $\sum_{i=k}^{m-1} a_i q_l h(y_{il}) = 0$ for scalars $a_i \in K$
with $a_k \ne 0$. Set $y = \sum_{i=k}^{m-1} a_i y_{il} \in \St(u_m)$.  First
we observe that $(k,l) \ne (1,0)$, for otherwise we would have $\rho h(y) =
a_k e + Je \ne 0$, which would make $h(y)$ a top element of
$A$, and hence guarantee that $q_l h(y) \ne 0$ by Proposition 16(A); but the
latter is inconsistent with our choice of $y$, and therefore we conclude
$(k,l) \ne (1,0)$.  This permits us to define an element $z \in \St(u_m)$ as
follows:  If $l = 0$, we have $k \ge 2$, which legitimizes the definition $z
=
\sum_{i=k}^{m-1} a_i y_{i-1,l}$; if, on the other hand, $l \ge 1$, we set $z
= \sum_{i=k}^{m-1} a_i y_{i,l-1}$.  We note that, in either case, the
nonzero scalar $a_k$ accompanies an element
$y_{ij}$ with $(i,j) < (k,l)$.  Hence, in case $l=0$, the minimal choice of
$(k,l)$ ensures that $p_{t+1} h(y) = q_t h(z) =
\sum_{i=k}^{m-1} a_i q_t h(y_{i-1,t}) \ne 0$, while, in case $l
\ge 1$, that choice yields $p_l h(y) = q_{l-1} h(z) =
\sum_{i=k}^{m-1} a_i q_{l-1} h(y_{i,l-1}) \ne 0$.  In other words,
$p_{t+1} h(y)$ is a nonzero element in $\soc A \cap q_t A$ in the first
case, and $p_l h(y)$ is a nonzero element in $\soc A \cap q_{l-1} A$ in the
second.  In either case, our hypothesis combines with part (B) of
Proposition 16 to show that $h(y)$ is a top element of type $e(l)$ of $A$,
whence our assumption that $q_l h(y)$ be zero contradicts part (A) of
Proposition 16.  We thus conclude that $\bigl( \la \etil(j) / J\etil(j)
\bigr)^{m-1}$ is contained in the socle of $A$ as claimed.

We have shown that, for each $j \in \{0, \dots, t\}$, all linear combinations
$\sum_{i=1}^{m-1} a_i q_j h(y_{ij})$ with $(a_1, \dots, a_{m-1}) \ne 0$ are
nonzero elements of the socle of $A$, and invoking again part (B) of
Proposition 16 on the model of the preceding paragraph, we infer that all
linear combinations
$\sum_{i=1}^{m-1} a_i  h(y_{ij})$ with $(a_1, \dots, a_{m-1}) \ne 0$  are top
elements of $A$.  This means that the elements
$h(y_{2k}),
\dots, h(y_{mk})$ are linearly independent modulo $JA$ and thus yields the
required containment of $\bigl( \la e(j) / Je(j) \bigr)^{m-1}$ in $A/JA$.

In conclusion, arbitrarily high finite powers of $\bigoplus_{i=0}^t
\la e(i)/Je(i)$ and $\bigoplus_{i=0}^t \la \etil(i)/J\etil(i)$ are
$\pinflamod$-phantoms, as well as $\sinflamod$-phantoms of $S$, and hence so
are their direct limits
$\bigl (\bigoplus_{i=0}^t \la e(i)/Je(i)\bigr)^{(\NN)}$ and
$\bigl( \bigoplus_{i=0}^t \la \etil(i)/J\etil(i) \bigr)^{(\NN)}$. \qed
\enddemo

Of course the mirror image of Lemma 23 relative to the central axis of $w$ is
also true, since replacing $w$ by its inverse is harmless. The argument we gave
for the lemma actually proves a little more than we stated in our conclusion: 
Namely, if
$u$ is a segment of $v$ such that $\St(u)$ has square-free socle (this is
for instance true for any syllable $u$ of $v$), then
$\bigl(\St(u)\bigr)^{(\NN)}$ is a $\pinflamod$-phantom of $S$. The part of
the lemma that addresses $\sinflamod$ will be superseded by the stronger
assertion of Theorem 5 that $\St(w)$ is an $\sinflamod$-phantom
of $S$.
\bigskip

\demo{\bf Proof of Theorem 5 via Proposition 14}  Let $e$ be one of the primitive
idempotents
$e_1, \dots, e_n$, and again denote the characteristic word of $S = \la e /
Je$ by $w  = (p_i^{-1}q_i)_{i \in \ZZ}$.

(I)  The word $w$ is centered at $e$ by construction.  For left and right
periodicity of $w$, we refer to the final paragraph of Construction 12.  To
justify the upper bound on the number of steps required to explicitly
determine $w$ from
$\Gamma$ and a set of paths generating $I$, we recall that  the principal
right segment of $w$ is completely determined, once we have constructed its
first $4n+1$ syllables $q_0, p_1^{-1}, q_1, p_2^{-1}, \dots, q_{2n}$. The
algorithm of \cite{\pre} for determining the projective dimensions of
$\la$-modules which are either uniserial or have graphs of the form `$V$',
moreover, allows us to find each successive pair $q_i p_i^{-1}$ of syllables
of
$w$ in $\le 2 (\dim_K \la)^3$ steps.
\medskip

(II) In light of Proposition 1, finiteness of the projective dimension of
$\St(w)$ is an immediate consequence of the construction.

To verify the second part of (II), let once more $\varphi: \St(w) \rightarrow S
=\la e / Je$ be the canonical map of the centered word
$w$.  In reference to the graph of $\St(w)$ at the beginning of Section 7,
we thus have $\varphi(x_i) = \delta_{i0} (e + Je)$ for
$i \in \ZZ$.   We first show that every homomorphism $M
\rightarrow S$, where $M$ is an object of $\sinflamod$, factors through
$\varphi$; this will prove effectiveness, once we have shown that $\St(w)$
is an $\sinflamod$-phantom of $S$.  It is clearly harmless to assume $M$ to
be indecomposable; in other words, we focus on the situation where $M$ is a
string module of finite projective dimension, based on a finite word $\what =
\hat{p}_0^{-1}\hat{q}_0 \dots \hat{p}_m^{-1}\hat{q}_m$ with $m \ge 0$, such
that the string module $M = \St(\what)$ has a standardized sequence $y_0,
y_1,\dots, y_m$ of $m+1$ top elements; see Section 2 for our conventions.
Then, clearly, those homomorphisms which send any given top element $y_i$ of
type $e$ to $e + Je$ and all the other $y_j$ to zero constitute a $K$-basis
of $\Hom_{\la}(\St(\what),S)$, whence we may assume that our map
$f \in \Hom_{\la}(\St(\what),S)$ is of the latter ilk; so suppose that there
exists $i \in \{1, \dots, m\}$ such that $y_i = e y_i$, and let $f \in
\Hom_{\la}(\St(\what),S)$ be as described.  We relabel the syllables of the
word $\what$ if necessary, to center it in the idempotent $e$ corresponding
to $y_i$, and to ensure that $\what$ has the same orientation as $w$. If the
standardized sequence of top elements $y_j$ is re-indexed accordingly, the
graph of $\St(\what)$ takes on the form

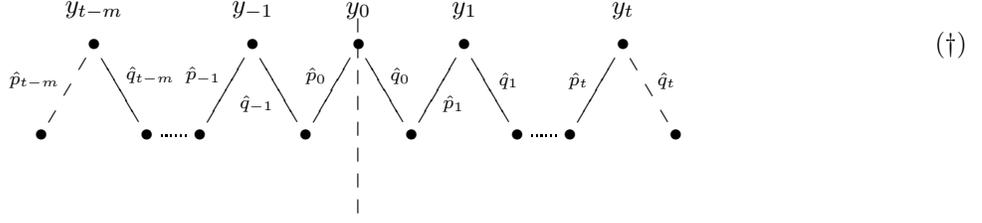
\ignore
$$\xymatrixrowsep{2pc}\xymatrixcolsep{0.75pc}
\xymatrix{
 &\bullet \dropup{3}{y_{t-m}} \ar@{--}[dl]_{\hatp_{t-m}}
\edge[dr]^(0.45){\hatq_{t-m}} &&&\bullet
\dropup{3}{y_{-1}} \edge[dl]_(0.45){\hatp_{-1}} \edge[dr]_(0.55){\hatq_{-1}}
&&\bullet
\save+<0ex,3ex> \drop{y_0} \ar@{--}[dd] \restore
\edge[dl]_(0.45){\hatp_0} \edge[dr]^(0.45){\hatq_0} &&\bullet \dropup{3}{y_1}
\edge[dl]^(0.55){\hatp_1} \edge[dr]^{\hatq_1} &&&\bullet \dropup{3}{y_t}
\edge[dl]_{\hatp_t}
\ar@{--}[dr]^{\hatq_t} \\
\bullet &&\bullet \ar@{.}[r] &\bullet &&\bullet &&\bullet &&\bullet
\ar@{.}[r] &\bullet &&\bullet\\
 &&&&&& } \tag\dagger$$
\endignore

\noindent for some $t \ge 0$. It is clearly innocuous to assume that
$\phat_t$ is nontrivial.

We now construct a map $g: \St(\what) \rightarrow \St(w)$ by starting with
the assignment $g(y_0) = x_0$.  Next, we inductively define the images
$g(y_j)$ for $1 \le j \le t$ in such a way that
$g(y_0), \dots, g(y_t)$ satisfy all of the relations tying the $\phat_j y_j$
and
$\qhat_j y_j$ for $0 \le j \le t$ together; and finally, we do the same for
$t-m \le j \le -1$.  Since the relations of $\St(\what)$ can be generated by
relations involving at most two consecutive $y_i$, this will ensure that our
assignments induce a well-defined homomorphism $g \in
\Hom_{\la}(\St(\what),St(w))$.

If none of the nontrivial paths $\qhat_i$,
$\phat_{i+1}$ for $i \ge 0$ is a right discontinuity of $\what$, the
assignments
$g(y_i) = x_i$ for $1 \le i \le t$ satisfy our requirements.  Next suppose
that the first right discontinuity of $\what$ is some nontrivial path
$\qhat_k$ with
$k \ge 0$.  Then the definitions $g(y_i) = x_i$ for $1 \le i \le k$ and
$g(y_i) = 0$ for $i > k$ are as required, for Lemma 20 tells us that
$q_k$ is a proper right subpath of $\qhat_k$, whence $\qhat_k g(y_k) = 0$.
If, finally, $\what$ has a first right discontinuity of the form $\phat_k$
for some $k > 0$, then $p_k =
 \phat_k \nu$ for a nontrivial path $\nu$, again by Lemma 20.  In that case
the assignments  $g(y_i) = x_i$ for
$i < k$, $g(y_k) = \nu x_k$, and $g(y_i) = 0$ for $i > k$ satisfy our
demands; indeed, if $\qhat_k$ is trivial, then $k = t$, and otherwise
$\qhat_k g(y_k) = 0$ in view of the nontriviality of $\nu$.

Assignments $g(y_j)$ for $t-m \le j \le -1$ which are compatible with the
pertinent relations are made symmetrically.  This procedure clearly leads to
a homomorphism $g$ with $\varphi g = f$.
\medskip

To show that $\St(w)$ is an $\sinflamod$-phantom of $S = \la e / Je$, we
consider, for each nonnegative integer $k$, the centered word $$u_k =
\tilde{p}_{-(k+1)}^{-1} q_{-(k+1)} p_{-k}^{-1} q_{-k}
\dots p_0^{-1} q_0 \dots p_k^{-1}q_k p_{k+1}^{-1}
\tilde{q}_{k+1};$$ here $\tilde{p}_{-(k+1)}= \frak{p}_{-(k+1)} p_{-(k+1)}$
and $\tilde{q}_{k+1} = \frak{q}_{k+1} q_{k+1}$ with
$\frak{p}_{-(k+1)}$ and $\frak{q}_{k+1}$ chosen as in Proposition 1, that
is, $\tilde{p}_{-(k+1)}$ and $\tilde{q}_{k+1}$ are the longest paths in
$K\Gamma \setminus I$ containing $p_{-(k+1)}$ and
$q_{k+1}$ as right subpaths.  By construction of $w = (p_i^{-1} q_i)_{i \in
\ZZ}$ and Proposition 1, the $u_k$ are words of finite projective
dimension.  As usual, we fix a standardized sequence of top elements
$y_{-(k+1)}, \dots, y_0, \dots, y_{k+1}$ of
$\St(u_k)$, and focus on the canonical map $\psi_k: \St(u_k)
\rightarrow S$, i.e., $\psi_k(y_i) = \delta_{i0} (e + Je)$.  We claim that
any module $A \in \sinflamod$ with the property that
$\psi_k$ factors through some map in $\Hom_{\la}(A,S)$ contains the string
module $\St(p_{-k}^{-1} q_{-k} \dots p_k^{-1} q_k)$ as a submodule. Once
established, this claim will entail that
$\St(p_{-k}^{-1} q_{-k} \dots p_k^{-1} q_k)$ is an
$\sinflamod$-phantom of $S$ for each $k$, whence so is the obvious direct
limit
$$\varinjlim_{k \in \NN} \St(p_{-k}^{-1} q_{-k} \dots p_k^{-1} q_k) =
\St(w).$$

Thus, in order to prove our claim, we suppose that $\psi_k$ factors through
$A \in \sinflamod$ and write $A = \bigoplus_{i=1}^r A_i$, where each
$A_i$ is a finite dimensional string module of finite projective dimension
and
$\pi_i: A \rightarrow A_i$ is the corresponding canonical projection.  That
$\psi_k$ factors through a map in $\Hom(A,S)$ clearly implies the existence
of a map $g \in \Hom_{\la}(\St(u_k), A)$ and an index $l$ such that
$\pi_l g(y_0)$ is a top element of $A_l$.  Suppose that $A_l = \St(\what)$,
where
$\what = \phat_0^{-1} \qhat_0 \dots \phat_m^{-1} \qhat_m$ is a word of
finite projective dimension such that $m \ge 0$ and $\St(\what)$ has a
standardized sequence of $m+1$ top elements, say $z_0, \dots, z_m$.

By \cite{\CraBoe, Theorem}, $\Hom_{\la}(\St(u_k),\St(\what))$ is generated,
as a $K$-space, by maps
$$h = h[\rho, \frak{u}_k,
\widehat{\frak{w}}]$$
 of the following ilk: $\frak{u}_k$ and
$\widehat{\frak{w}}$ are subgraphs of the graphs of the words
$u_k$ and $\what$, respectively (see Section 2 for our conventions), the
first closed under arrows whose endpoints belong to $\frak{u}_k$, the second
closed under arrows whose starting points belong to $\widehat{\frak{w}}$;
moreover, $\rho$ denotes an isomorphism $\frak{u}_k \rightarrow
\widehat{\frak{w}}$ of directed graphs, sending any arrow in $\frak{u}_k$ to
an arrow in
$\widehat{\frak{w}}$ that carries the same label, such that the homomorphism
$h$ is induced by $\rho$ (in the only meaningful way).  Note that the image
of $h$ is $\St(\widehat{\frak{w}})$, the latter being a submodule of
$\St(\what)$ by the closure condition imposed on $\widehat{\frak{w}}$.   So
the fact that there exists a map in $\Hom_{\la}(\St(u_k),\St(\what))$ sending
$y_0$ to a top element of $\St(\what)$ ensures the existence of a triple
$[\rho, \frak{u}_k, \widehat{\frak{w}}]$, as described, together with an
index $i \in\{0,\dots, m\}$, satisfying the following requirements: 
$\widehat{\frak{w}}$ includes that vertex in the graph of $\what$ which
corresponds to the top element $z_i$ of $\St(\what)$  --  call that vertex
$z$ (it is the joint starting vertex of the paths $\phat_i$ and $\qhat_i$ in
the graph of $\what$); moreover, $\frak{u}_k$ includes that vertex in the
graph of $u_k$ which corresponds to the top element $y_0$ of
$\St(u_k)$ (namely, the joint starting vertex of $p_0$ and $q_0$ in the
graph of $u_k$); and, finally, $\rho$ sends this latter vertex to $z$.

In light of the preceding discussion, it suffices to prove that
$\widehat{\frak{w}}$ contains a subgraph isomorphic to that of the word
$p_{-k}^{-1} q_{_k} \dots p_k^{-1} q_k$.   As in the effectiveness proof
above, we adjust the labeling of the syllables of $\what$ and the
standardized sequence of top elements of
$\St(\what)$ so that $\what$ becomes a word which is centered at
$i=0$, and the graph of $\St(\what)$ has the form ($\dagger$), displayed at
the outset of our proof of (II), for a suitable nonnegative integer $t$.
Moreover, it is clearly harmless to assume that $\what$ has the same
orientation as $w$. Then
$\phat_0$ contains $p_0$ as a right subpath, and $\qhat_0$ contains $q_0$ as
a right subpath by Proposition 16(A), because
$\what$ is a word of finite projective dimension. Since the graph
$\widehat{\frak{w}}$ contains the paths $\phat_0$ and $\qhat_0$ due to its
closure property, we deduce that $\phat_0 = p_0$ and
$\qhat_0 = q_0$. Part (B) of that same proposition now tells us that the
paths $\phat_{-1}$ and $\phat_1$ are contained in
$p_{-1}$ and $p_1$ as left subpaths, respectively.  If $p_{-1}$ and $p_1$
are trivial, we are done.  So we assume that $p_1$ is nontrivial and write
$p_1 = \phat_1 \nu$, where $\nu$ is a path of length $\ge 0$.  Then $\qhat_0
z_0 = q_0 z_0 = \phat_1 \nu z_1$ is a nonzero element of $\soc \St(\what)$,
whence Proposition 16(B) guarantees that $\nu z_1$ is a top element of
$\St(\what)$.  We infer that $\nu$ is trivial and that the vertex of $\what$
corresponding to $z_1$ also belongs to $\widehat{\frak{w}}$. Consequently,
the preceding argument can be duplicated to show
$\qhat_1 = q_1$, and next, in case $p_2$ is nontrivial, the equalities
$\phat_2 = p_2$ and $\qhat_2 = q_2$. An obvious induction on $i \in \{1,
\dots, k\}$ yields $\phat_i = p_i$ and
$\qhat_i = q_i$, whenever $p_i$ is nontrivial, and an analogous argument
applies to the paths $\phat_i$ and $\qhat_i$ with negative indices $i \in
\{-k, \dots, -1\}$.  This shows that the graph of the word $p_{-k}^{-1}
q_{-k} \dots p_k^{-1} q_k$ is indeed a subgraph of $\widehat{\frak{w}}$,
thus proving our claim and finishing the proof of part (II) of Theorem 5.
\medskip

(III)  The first assertion is an immediate consequence of the equivalences.
To prove the equivalences, we fix $k$ and write $S_k = S = \la e / Je$.
Moreover, we let $w_k = w = (p_i^{-1}q_i)_{i \in \ZZ}$ be the characteristic
word of
$S$ and $\varphi: \St(w) \rightarrow S$ the canonical homomorphism defined by
$\varphi(x_i) = \delta_{i0} (e + Je)$.

`(i)$\implies$(iii)'.  Suppose that $w$ is finite.  We need to ascertain
that every homomorphism $f \in \Hom_{\la}(M,S)$, where $M$ is an
indecomposable object of $\pinflamod$, factors through $\varphi$.  Once we
know that
$\varphi$ is a $\pinflamod$-approximation of
$S$, minimality will be automatic, because
$\St(w)$ is known to be indecomposable (see \cite{\BuRi}). In case $M$ is a
string module of finite projective dimension, the required factorization
property follows from the fact that $\varphi$ is an {\it effective}
$\sinflamod$-phantom of $S$, which was established in part (II).

Now we focus on a band module $M = \Bd(\vhat^r, \phi)$ of finite projective
dimension, where $\vhat =
\hat{p}_0^{-1}\hat{q}_0 \dots \hat{p}_t^{-1}\hat{q}_t$ is a primitive word
with all of the listed syllables nontrivial, and again denote by
$y_{10}, \dots, y_{1t}$, $y_{20}, \dots, y_{2t}$, $y_{r0}, \dots, y_{rt}$ a
standardized sequence of top elements of $M$.  It clearly suffices to show
that any map $f \in
\Hom_{\la}(M,S)$ which sends precisely one of the top elements
$y_{ij}$ of type $e$ to $e + Je$ and the others to zero factors through
$\varphi$ (if none of the $y_{ij}$ is of type $e$, our requirement is
void).  So let us assume that $f : M
\rightarrow S$ is a map of the described ilk. Moreover, it is harmless to
adjust the setup so that $f(y_{10}) = e +Je$ and
$f(y_{ij}) = 0$ for $(i,j) \ne (1,0)$.  In this situation, $M$ has the same
graph as the pseudo-band module $B$ in the statement of Lemma 22.

As in Lemmas 20 and 22, we denote by $\what$ the twosided infinite word
$\dots \vhat \vhat \vhat \dots$, which we again assume to be centered at $e$
and have the same orientation as $w$.  Clearly,
$\what$ has right and left discontinuities, since the characteristic word
$w$ of $S$ is finite by hypothesis, whereas
$\what$ is twosided infinite.  In light of Proposition 1, moreover, $\what$
is a word of finite projective dimension, and so Lemma 22 applies to
guarantee that $f$ factors through $\varphi$ as required.

The implications `(iii)$\implies$(v)$\implies$(iv)' and
`(iii)$\implies$(ii)' are obvious.

`(iv)$\implies$(i)' follows from the fact that $\St(w)$ is an
$\sinflamod$-phantom of $S$ by part (II).

`(ii)$\implies$(i)'.  Suppose  $w = (p_i^{-1} q_i)_{i \in \ZZ}$ is
infinite;  w\.l\.o\.g\., the principal right segment of $w$ is infinite and
so, in particular, the syllable $q_0$ is nontrivial.

First assume that there exists an index $j > 0$ such that $q_j = q_0$, and
let $t \ge 1$ be minimal with  the property that $q_t = q_0$.  Then the
principal right segment $w_{\text right}$ of $w$ has the form
$$w_{\text right}  = (q_0 p_1^{-1} \dots q_t p_{t+1}^{-1}) (q_0 p_1^{-1}
\dots q_t p_{t+1}^{-1}) (q_0 p_1^{-1} \dots q_t p_{t+1}^{-1})
\dots$$ by Construction 12; in other words, if $v = p_{t+1}^{-1} q_0 \dots
p_t^{-1} q_t$ is a primitive word such that
$p_{t+1}^{-1} w_{\text right} = \dots vvv \dots$,  Lemma 23 tells us that
$S$ fails to have a $\pinflamod$-approximation. In case
$p_0$ is nontrivial and there exists an index $j < 0$ with $p_j = p_0$, the
situation is symmetric  --  just flip the characteristic word $w$ about its
central axis.

Now suppose that, for all $j > 0$, we have $q_j \ne q_0$, and $p_0$ is
either trivial, or else $p_j \ne p_0$ for all $j < 0$.  Assume that, to the
contrary of our claim,
$S$ has a $\pinflamod$-approximation.  Our aim is to infer that then $S$ has
an
$\sinflamod$-approximation as well; but this is incompatible with the
already established implication `(iv)$\implies$(i)'.  To that end, let
$A$ be any $\pinflamod$-approximation of $S$ and
$B$ a band module occurring as a direct summand of $A$, say $A = B \oplus
C$.  We will construct another
$\pinflamod$-approximation of $S$ which has the form $\bigoplus_{\text
finite} \St(u_j) \oplus C$, for suitable words $u_j$ of finite projective
dimension.  In this way, we can eliminate all band module summands from $A$
in favor of direct sums of string modules, to arrive at another
$\pinflamod$-approximation of $S$, this one belonging to
$\sinflamod$.  This will then give us the desired contradiction.

In order to replace the summand $B$ of $A$ by a direct sum of string modules
as indicated, write $B = Bd(\vhat ^r, \phi)$, where
$\vhat = \phat_0^{-1} \qhat_0 \dots \phat_t^{-1} \qhat_t$ is a primitive
word, $r \ge 1$. For a graph of $B$, we refer to the graph of the
pseudo-band module of the same name in the statement of Lemma 22; again, we
let $y_{10}, \dots, y_{1t}, \dots, y_{r0},
\dots, y_{rt}$ be the corresponding standardized sequence of top elements of
$B$.  As in the argument for `(i)$\implies$(iii)', it suffices to prove the
following:  Given any top element $y \in
\{y_{10}, \dots, y_{rt}\}$ of type $e$ of $B$, the map $f: B
\rightarrow S$, sending $y$ to the residue class $e + Je$ and the other top
elements $y_{ij}$ to zero, factors through a map $\St(u)
\rightarrow S$, where $u$ is a suitable word of finite projective
dimension.  Moreover, it is harmless to assume that $y$ is the top element
corresponding to the joint starting vertex of the paths
$\phat$ and $\qhat$. Once more, we center the infinite periodic word $\what
= \dots \vhat \vhat \vhat \dots$ in some occurrence of this vertex, and
assume (clearly still without losing generality) that $\what$ has the same
orientation as $w$.  Since $B$ has finite projective dimension, so does the
word $\what$, by Proposition 1. The non-periodicity conditions we imposed on
$w$, moreover, force left and right discontinuities on $\what$.  In light of
Lemma 22, our test map $f$ thus factors through the canonical homomorphism
$\varphi: \St(w) \rightarrow S$.  Say $f =
\varphi g$, for a suitable map $g \in \Hom_{\la}(B, \St(w))$.  Let
$w' = (p_i^{-1} q_i)_{-m \le i \le m}$ be a finite segment of $w$ such that
$\St(w')$ is a submodule of $\St(w)$ containing the image of $g$. Moreover,
let $u$ be a finite word of finite projective dimension which in turn
contains $w'$ as a segment; such a word $u$ exists by Observations 10.  Then
$f = \varphi h$, where $h: B \rightarrow \St(u)$ denotes the map resulting
from $g$ through restriction of the range, which finishes the argument we
have laid out.

\smallskip To prove the supplementary statement, let $w$ be finite. Suppose
that the $K$-dimension of $\St(w)/J\St(w)$ exceeds $4n$, where $n$ is the
number of vertices of $\Gamma$.  Then either $p_{2n}$ or
$q_{-2n}$ is nontrivial.  We may assume the former path to be nontrivial,
the other case leading to a symmetric situation.   The string module
$\St(q_0 p_1^{-1}\dots q_{2n-1} p_{2n}^{-1})$ then has $2n + 1$ top elements
$x_i$, which are $K$-linearly independent modulo the radical.  At least
three of these are normed by the same primitive idempotent.  Hence at least
two of the latter, say $x_k$ and $x_l$ for suitable $k < l$, reside atop
paths $p_k$, $p_l$ starting in the same arrow. In this situation,
Construction 12 yields $q_k = q_l$, $p_{k+1} = p_{l+1}$, etc. (consult
Observations 11), and consequently all the $p_i$ with positive index are
nontrivial.  But this makes the word $w$ infinite, thus contradicting our
hypothesis.

The proof of Theorem 5 is thus complete. \qed
\enddemo

\head  8. Concluding remarks \endhead

Theorem 5 and Proposition 14 tell us that, if $w = w(S)$ is the
characteristic word of the simple module $S \in \lamod$, the string module
$\St(w)$ is an effective $\sinflamod$-phantom of $S$.  On the other hand,
$\St(w)$ is not an effective
$\pinflamod$-phantom  of $S$ in general.  The first known example of a
finite dimensional algebra $\la$ for which $\pinflamod$ fails to be
contravariantly finite, presented in \cite{\IgSmTo}, already demonstrates
this.

\example{Example 23}  Let $\la$ be the string algebra whose two
indecomposable projective left modules have the following graphs:

\ignore
$$\xymatrixrowsep{1pc}\xymatrixcolsep{0.5pc}
\xymatrix{
 &1 \edge[dl]_\alpha \edge[dr]^\beta &&&2 \edge[d]\\ 2 \edge[d] &&2 &&1\\ 1
}$$
\endignore

\noindent Using the techniques we developed, we see that the characteristic
phantom $\St(w_1)$ of $S_1$ has graph

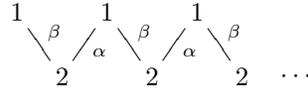
\ignore
$$\xymatrixrowsep{1pc}\xymatrixcolsep{0.5pc}
\xymatrix{ 1 \edge[dr]^\beta &&1 \edge[dl]^\alpha \edge[dr]^\beta &&1
\edge[dl]^\alpha \edge[dr]^\beta\\
 &2 &&2 &&2 &\cdots }$$
\endignore

\noindent  This infinite dimensional module is known to be both an
$\sinflamod$- and a
$\pinflamod$-phantom of $S_1$ (see \cite{\HaHZ}), but, while it is effective
as an $\sinflamod$-phantom by Theorem 5, it is not effective as a
$\pinflamod$-phantom.  Indeed, for the band module $M = \la e_1 /
\la(\alpha - \beta) \in \pinflamod$, the canonical epimorphism $f: M
\rightarrow S_1$ clearly fails to factor through $\St(w_1)$.

Furthermore, we note that the string module $\St(w_1)$, unique relative to the
given coordinatization of $\la$, does depend on the coordinates. If we adopt the
new arrows $\alpha'= \alpha - \beta$ and $\beta'=
\beta$ from the vertex $1$ to the vertex $2$, the characteristic phantom
$\St(w')$ of $S_1$ with respect to the new coordinatization of $\la$ has
graph

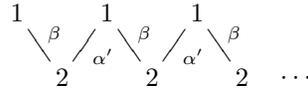
\ignore
$$\xymatrixrowsep{1pc}\xymatrixcolsep{0.5pc}
\xymatrix{ 1 \edge[dr]^\beta &&1 \edge[dl]^{\alpha'} \edge[dr]^\beta &&1
\edge[dl]^{\alpha'} \edge[dr]^\beta\\
 &2 &&2 &&2 &\cdots }$$
\endignore

\noindent  It is readily seen that $\St(w_1)$ is  not isomorphic to
$\St(w_1')$.  \qed
\endexample

In Part II, we will address the problem of supplementing the effective
$\sinflamod$-phantoms $\varphi_i: \St(w_i) \rightarrow S_i$ to effective
$\pinflamod$-phantoms.  It will turn out that all information required to
construct such `completions' of the
$\St(w_i)$ is stored in the characteristic words $w_i$. 

Concerning the shape of the minimal $\pinflamod$-approximations of the
$S_i$ in case the word $w = w(S)$ is finite, Theorem 5 guarantees
that the top of the minimal
$\pinflamod$-approximation $\St(w)$ of $S$ has $K$-dimension at most $4n$,
where $n$ is the number of distinct simple
$\la$-modules.  This bound stems from the fact that the multiplicity of any
simple module $T$ in $\St(w) / J \St(w)$ is bounded above by $4$ (see
Construction 12).  While the bound on
$\dim_K \bigl(\St(w) / J \St(w) \bigr)$ can be tightened with some
additional effort, the bound on multiplicities in $\St(w) / J \St(w)$ is sharp.

\example{Example 24}  Let $\la$ be the string algebra whose indecomposable
projective left modules have the following graphs

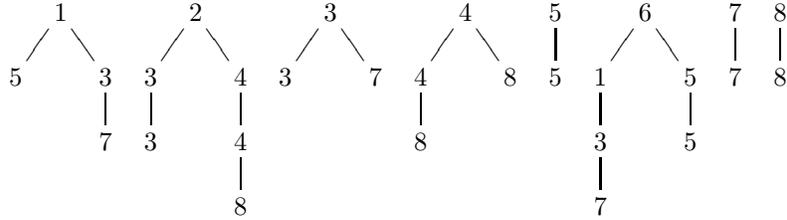
\ignore
$$\xymatrixrowsep{1pc}\xymatrixcolsep{0.5pc}
\xymatrix{
 &1 \edge[dl] \edge[dr] &&&2 \edge[dl] \edge[dr] &&&3 \edge[dl] \edge[dr]
&&&4 \edge[dl] \edge[dr] &&5 \edge[d] &&6 \edge[dl] \edge[dr] &&7
\edge[d] &8 \edge[d]\\ 5 &&3 \edge[d] &3 \edge[d] &&4 \edge[d] &3 &&7 &4
\edge[d] &&8 &5 &1
\edge[d] &&5 \edge[d] &7 &8\\
 &&7 &3 &&4 \edge[d] &&&&8 &&&&3 \edge[d] &&5\\
 &&&&&8 &&&& &&&&7 }$$
\endignore

\noindent  Then $\pinflamod$ is contravariantly finite in $\lamod$, as can
readily be ascertained with the aid of Theorem 5, and the graph of the
minimal
$\pinflamod$-approximation of
$S_1$ is

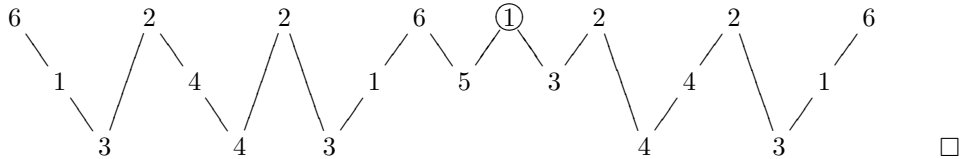
\ignore
$$\xymatrixrowsep{1pc}\xymatrixcolsep{0.5pc}
\xymatrix{ 6 \edge[dr] &&&2 \edge[ddl] \edge[dr] &&&2 \edge[ddl] \edge[ddr]
&&&6
\edge[dl] \edge[dr] &&\encirc{1} \edge[dl] \edge[dr] &&2 \edge[dl]
\edge[ddr] &&&2 \edge[dl] \edge[ddr] &&&6 \edge[dl]\\
 &1 \edge[dr] &&&4 \edge[dr] &&&&1 \edge[dl] &&5 &&3 &&&4 \edge[dl] &&&1
\edge[dl]\\
 &&3 &&&4 &&3 &&&&&&&4 &&&3 &&&&\square }$$
\endignore
\endexample

This final example also exhibits the necessity of recording the center of
the characteristic word $w$ of $S$, if one aims at pinning down the correct
homomorphism $\varphi: \St(w) \rightarrow S$, through which all
homomorphisms in $\Hom_{\la} \bigl(\sinflamod,S \bigr)$ will factor. 
Indeed, the minimal
$\pinflamod$-approximation of $S_2$ is the string module with
graph

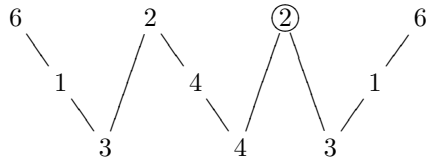
\ignore
$$\xymatrixrowsep{1pc}\xymatrixcolsep{0.5pc}
\xymatrix{ 6 \edge[dr] &&&2 \edge[ddl] \edge[dr] &&&\encirc{2} \edge[ddl]
\edge[ddr] &&&6
\edge[dl] \\
 &1 \edge[dr] &&&4 \edge[dr] &&&&1 \edge[dl] \\
 &&3 &&&4 &&3 }$$
\endignore

\noindent where again the center is highlighted.

\enddocument